\documentclass[graybox]{svmult}

\usepackage{mathptmx}       
\usepackage{helvet}         
\usepackage{courier}        
\usepackage{type1cm}        
\usepackage{makeidx}         
\usepackage{graphicx}        
\usepackage{multicol}        
\usepackage[bottom]{footmisc}

\makeindex             

\usepackage[T1]{fontenc}

\newif\iffinal

\finaltrue

\usepackage{siunitx}
\DeclareSIUnit\px{px}

\usepackage{dsfont}

\usepackage{amsmath}
\usepackage{amssymb}

\DeclareMathOperator{\dist}{dist}
\DeclareMathOperator{\pers}{pers}

\DeclareMathOperator{\persG}{\ensuremath{\pers_G}}

\newcommand{\diagram}  {\ensuremath{\mathcal{D}}}
\newcommand{\graph}    {\ensuremath{\mathcal{G}}}

\newcommand{\pixels}   {\ensuremath{\mathcal{P}}}

\usepackage{tikz}
\usepackage{pgfplots}

\usetikzlibrary{%
  arrows,
  chains,
  external,
  shapes,
}

\tikzsetexternalprefix{Figures/External/}

\pgfplotsset{compat=1.14}

\definecolor{created}   {RGB}{144,  0, 32}
\definecolor{destroyed} {RGB}{ 70,130,180}
\definecolor{irregular} {RGB}{255,191,  0}
\definecolor{persisting}{RGB}{150,150,150}

\definecolor{cardinal}{RGB}{196, 30, 58}

\DeclareRobustCommand{\inlinedot}[1]{\raisebox{0.25\height}{\tikz\fill[#1] (0,0) circle (0.50ex);}}

\usepackage{subfigure}

\newcommand{\subfigureCaptionSkip}{\vspace{-10pt}}%

\usepackage{paralist}
\usepackage{xspace}

\newcommand  {\nd}{\textsuperscript{\textup{nd}}\xspace}

\newcommand{\comment}[1]{}
\usepackage[absolute]{textpos}

\begin{document}

\title*{Persistence Concepts for 2D Skeleton Evolution Analysis}

\author{Bastian Rieck \and Filip Sadlo \and Heike Leitte}

\institute{Bastian Rieck \and Heike Leitte \at TU Kaiserslautern, \email{\{rieck, leitte\}@cs.uni-kl.de}
  \and Filip Sadlo \at Heidelberg University, \email{sadlo@uni-heidelberg.de}%
}

\maketitle

\abstract{%
In this work, we present concepts for the analysis of the evolution of
two-dimensional skeletons. By introducing novel persistence concepts, we
are able to reduce typical temporal incoherence, and provide insight
in skeleton dynamics. We exemplify our approach by means of a simulation of viscous
fingering---a highly dynamic process whose analysis is a
hot topic in porous media research.
}

\iffinal
\else
  \begin{textblock*}{\paperwidth}[1, 0](\paperwidth,0.5cm)
  \scriptsize%
  Authors' copy. Please refer to \emph{Topological Methods in Data Analysis and Visualization
  V: Theory, Algorithms, and Applications} for the definitive version of this chapter.
  \end{textblock*}
\fi

\section{Introduction}

\begin{figure}[tbp]
  \centering
  \subfigure[$t=8$\label{fig:viscfing_1-t8}]{\includegraphics[width=0.325\linewidth]{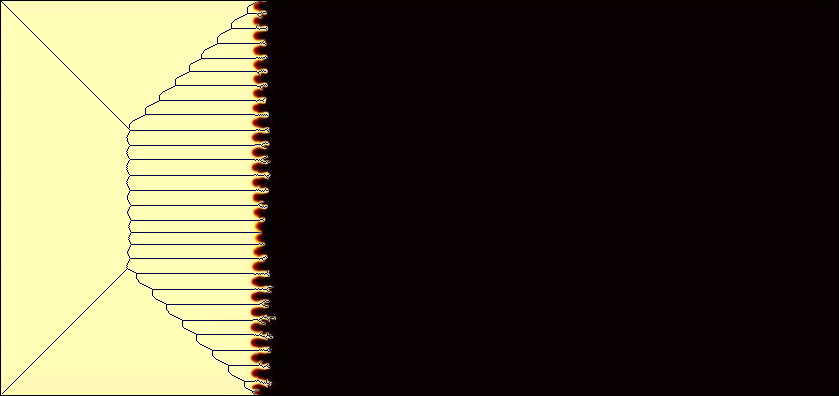}}\hfill%
  \subfigure[$t=30$\label{fig:viscfing_1-t30}]{\includegraphics[width=0.325\linewidth]{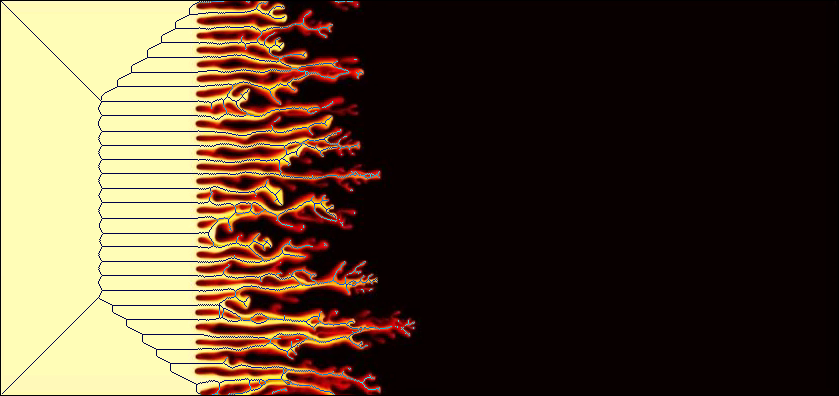}}\hfill%
  \subfigure[$t=70$\label{fig:viscfing_1-t70}]{\includegraphics[width=0.325\linewidth]{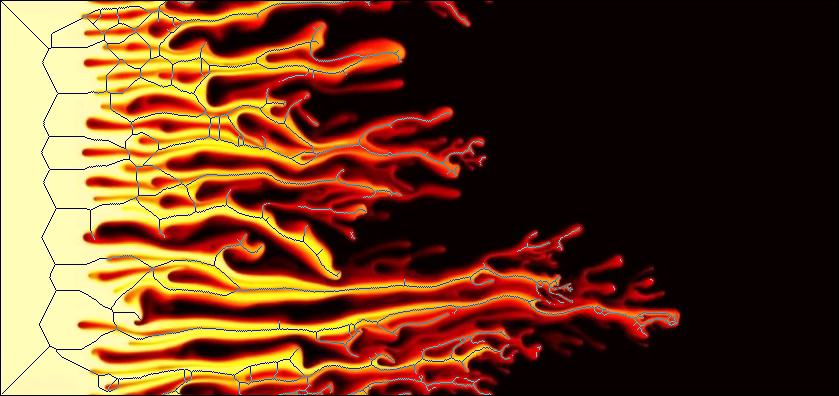}}%
  \subfigureCaptionSkip
  \caption{%
    Selected time steps~\subref{fig:viscfing_1-t8}--\subref{fig:viscfing_1-t70} of a 2D viscous
    fingering simulation~\cite{viscousFingeringVideo}, with extracted skeleton~(overlay).
    We used a conservative threshold for segmentation to suppress dark-red parts.
  }
  \label{fig:Skeletons overlay}
\end{figure}

There are many research problems that express themselves more in terms
of topological structure than morphology. Typical examples of such
processes include electrical discharge, the growth of crystals, and
signal transport in networks.
In this paper, we address \emph{viscous fingering}, where the interface between two
fluids is unstable and develops highly-dynamic ``finger-like'' structures.
A prominent cause for such structures are setups where a fluid with lower
viscosity~(Fig.~\ref{fig:viscfing_1-t8}--\subref{fig:viscfing_1-t70}, left) is injected into a fluid with higher
viscosity~(Fig.~\ref{fig:viscfing_1-t8}--\subref{fig:viscfing_1-t70}, right).
To analyze these processes, a straightforward approach employs traditional skeletonization
techniques for extracting the topology of each time step independently. Here, we employ iterative
thinning~\cite{thinning1992}. However, like all skeletonization techniques, the resulting skeletons
tend to be to temporally incoherent because the extraction is susceptible to small variations and
noise.
We present persistence concepts to address these issues and provide insight into the underlying
processes.

\section{Viscous Fingering}
%
Even though the methods described in this paper are generically applicable to time-varying
skeletons, we focus our analysis on skeletons that we extracted from \emph{viscous fingering}
processes. The term viscous fingering refers to the formation of structural patterns that appear
when liquids of different viscosity are mixed.
Under the right conditions, e.g., when water is being injected into glycerine, branch-like
structures---the eponymous \emph{viscous fingers}---begin to appear and permeate through the liquid of
higher viscosity.
Understanding the formation of these patterns is a prerequisite for the description of many natural
processes, such as groundwater flows. Consequently, researchers are interested in setting up
simulations that closely match the observations of their experiments.

Since each simulation uses a different set of parameters, summary statistics and comparative
visualizations are required in order to assess how well a simulation describes an experiment.
As a first step towards analyzing these highly-complex dynamics, we extract skeletons for each time
step of a simulation or an experiment. In this paper, we introduce several concepts for assessing
the inherent dynamics of these skeletons, permitting a comparative analysis.

\runinhead{Other methods}
%
In the context of analyzing viscous fingering, several other techniques exist. An approach by
Lukasczyk et al.~\cite{Lukasczyk17}, for example, uses tracking graphs to visualize the
spatio-temporal behavior of such processes. In a more general context, discrete Morse theory could
be applied to detect persistent structures in gray-scale images~\cite{Delgado-Friedrichs15}.
The applicability of these approaches hinges on the data quality, however. Our experimental data
suffers from a high noise level in which many smaller fingers cannot be easily identified by the
other approaches. This is why we decided to focus on conceptually simpler skeletonization techniques
for now.

\begin{figure}[tbp]
  \centering
  \iffinal
    \includegraphics{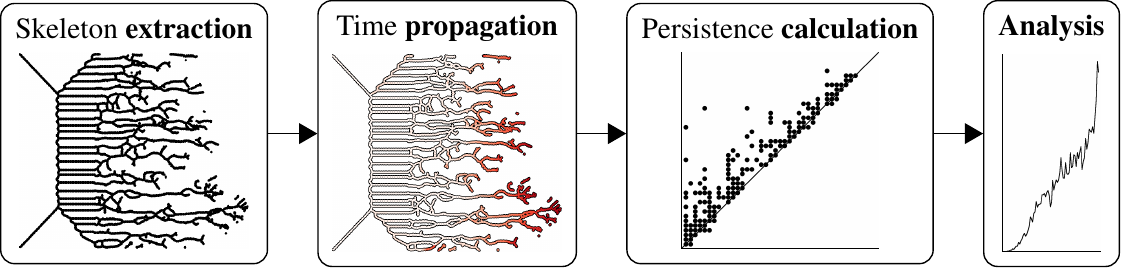}
  \else
    \begin{tikzpicture}[%
      start chain  =going right,
      node distance=5mm,]
      \tikzset{%
        block/.style = {%
          draw,
          on chain,
          on grid,
          rounded corners,
          align          = center,
          inner sep      = 0.5em,
          shape          = rectangle,
          minimum height = 3em,
          minimum width  = 1cm,
        },
        line/.style = {%
          >=triangle 60,
          draw,
          ->,
        }
      }
      \node[block] (Extraction) {%
        Skeleton \textbf{extraction}\\[0.1cm]%
        \includegraphics[height=2cm]{Figures/Skeleton_glyph}
      };

      \node[block] (Propagation) {%
        Time \textbf{propagation}\\[0.1cm]%
        \includegraphics[height=2cm]{Figures/Skeleton_ages_glyph}
      };

      \node[block] (Calculation) {%
        Persistence \textbf{calculation}\\[0.1cm]%
        \includegraphics[height=2cm]{Figures/Persistence_diagram_glyph}
      };

      \node[block] (Analysis) {%
        \textbf{Analysis}\\[0.1cm]%
        \includegraphics[height=2cm]{Figures/Analysis_glyph}
      };

      \path[line] (Extraction)  -- (Propagation);
      \path[line] (Propagation) -- (Calculation);
      \path[line] (Calculation) -- (Analysis);
    \end{tikzpicture}
  \fi
  \caption{%
    The basic pipeline of our approach. The first step, i.e., skeleton extraction, strongly depends
    on the desired application.
    Likewise, the analysis step can comprise different diagrams, summary statistics, and goals.
    Individual parts of the pipeline are replaceable, making our approach highly generic.
    Our current implementation uses an algorithm by Zhang and Suen~\cite{Zhang84} for skeleton
    extraction~(Sec.~\ref{sec:Propagation}).
    The subsequent propagation of creation times between time steps along all branches of the
    skeleton uses the methods described in the same section.
    From this extended skeleton, Sec.~\ref{sec:Persistence concepts} describes how to derive numerous
    persistence diagrams.
    Following this, we define multiple activity indicators based on these diagrams in Sec.~\ref{sec:Activity
    indicators}. Finally, Sec.~\ref{sec:Analysis} presents an analysis of different data sets under different aspects.
  }
  \label{fig:Pipeline}
\end{figure}

\section{Overview and Methods}

In this paper, we implement a pipeline that comprises the whole range of the analysis
process of a series of time-varying skeletons.
Fig.~\ref{fig:Pipeline} shows a schematic illustration and points to the corresponding
sections in which individual parts are described.
We provide an open-source implementation~(in \texttt{Python}) of the pipeline on
GitHub\footnote{\url{https://github.com/Submanifold/Skeleton\_Persistence}}. The repository includes
all examples, data, and instructions on how to reproduce our experiments.
For the analysis of our persistence diagrams, we implemented tools that build upon
\texttt{Aleph}\footnote{\url{https://github.com/Submanifold/Aleph}}, an open-source library for
persistent homology calculations.
We stress that our implementation is a proof of concept. Its computational bottleneck is the
brute-force matching~(which could be improved by using an approximate matching algorithm) that is
required as a precursor to creation time propagation.
More precisely, calculating all matches over all time steps takes between \SI{2}{\hour} and
\SI{6}{\hour}, while the subsequent propagation of creation times takes \SI{82}{\second}~(example
data, $\SI{839}{\px}\times\SI{396}{\px}$, 84 time steps), \SI{384}{\second}~(measured data,
$\SI{722}{\px}\times\SI{1304}{\px}$, 58 time steps), and \SI{524}{\second}~(simulation data,
$\SI{1500}{\px}\times\SI{1000}{\px}$, 37 time steps).
Finally, persistence diagram creation requires \SI{100}{\second}~(example data),
\SI{183}{\second}~(simulation data), and \SI{926}{\second}~(measured data), respectively.
The time for calculating activity indicators~(Sec.~\ref{sec:Activity indicators}), e.g., total
persistence, is negligible, as the persistence diagrams only contain a few hundred points. Please refer to Sec.~\ref{sec:Analysis} for more information about the individual data sets.

Subsequently, we will first briefly discuss skeleton extraction---both in terms of sets of pixels as well
as in terms of graphs.
Next, we explain the necessary steps for obtaining information about the ``creation time'' of pixels
and how to propagate said information over all time steps in order to obtain evolution information.
Based on this, we derive and exemplify several concepts motivated by topological persistence.

\subsection{Skeleton Extraction and Propagation of Pixel Creation Time}
\label{sec:Extraction and propagation}

Iterative thinning provides skeletons from binary images in a pixel-based format.
A sequence of skeletons thus gives rise to a sequence of pixel sets $\pixels_0$, $\pixels_1$, \dots,
$\pixels_k$, each corresponding to a time step $t_0$, $t_1$, \dots, $t_k$.
We employ \mbox{8-neighborhood} connectivity around each pixel, i.e., the set of all neighbors
including the diagonal ones, to convert each pixel set $\pixels_i$ into a graph~$\graph_i$.
Depending on the degree~$d$ of each vertex in~$\graph_i$, we can classify each pixel as being either
a regular point~($d=2$), a start/end point~($d=1$), or a branch point~($d\geq{}3$).  This also
permits us to define \emph{segments} formed by connected subsets of regular pixels.

\subsubsection{Pixel Matching}
\begin{figure}[bt]
  \centering
  \subfigure[$t=72$]{%
    \includegraphics[height=3.75cm]{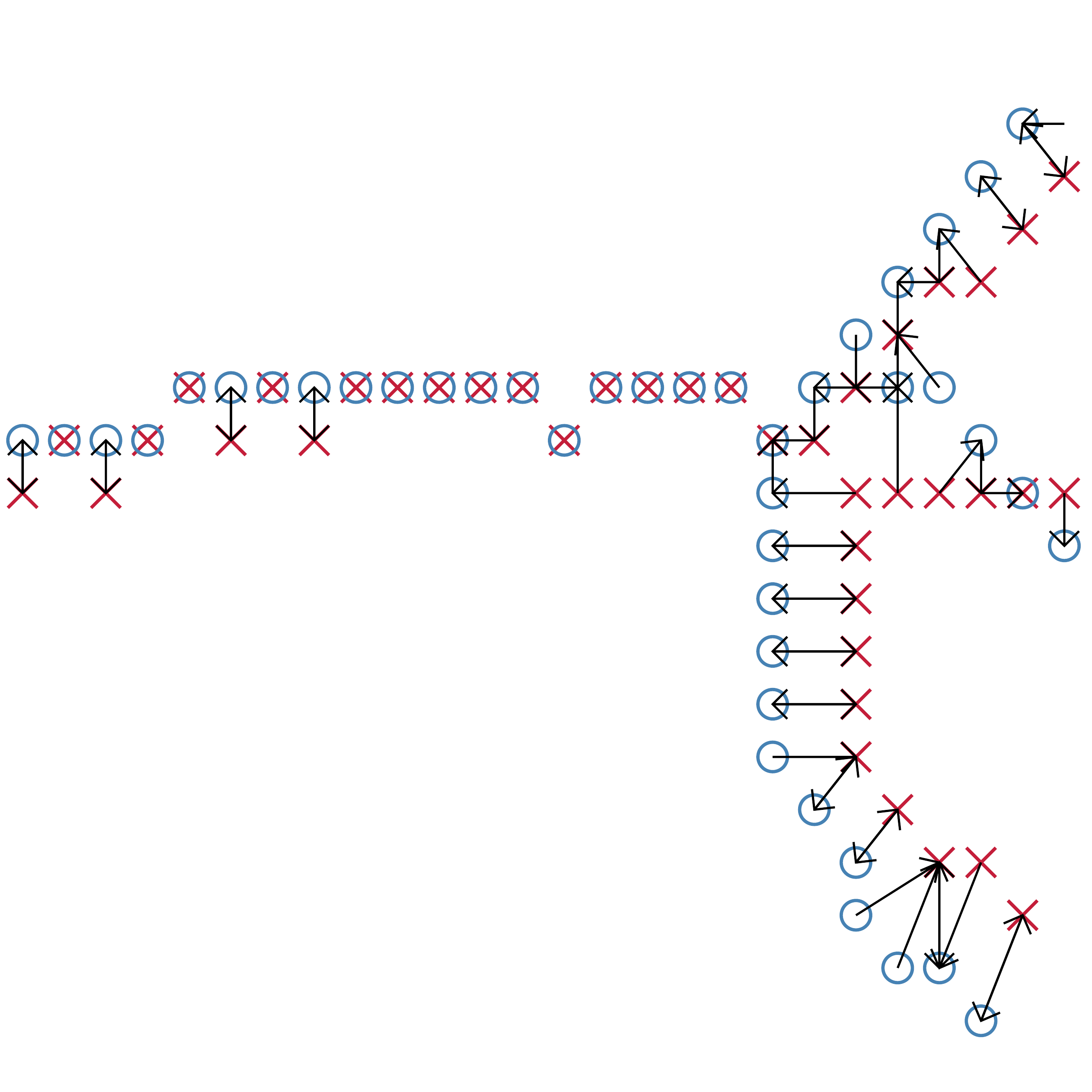}
  }
  \subfigure[$t=73$]{%
    \includegraphics[height=3.75cm]{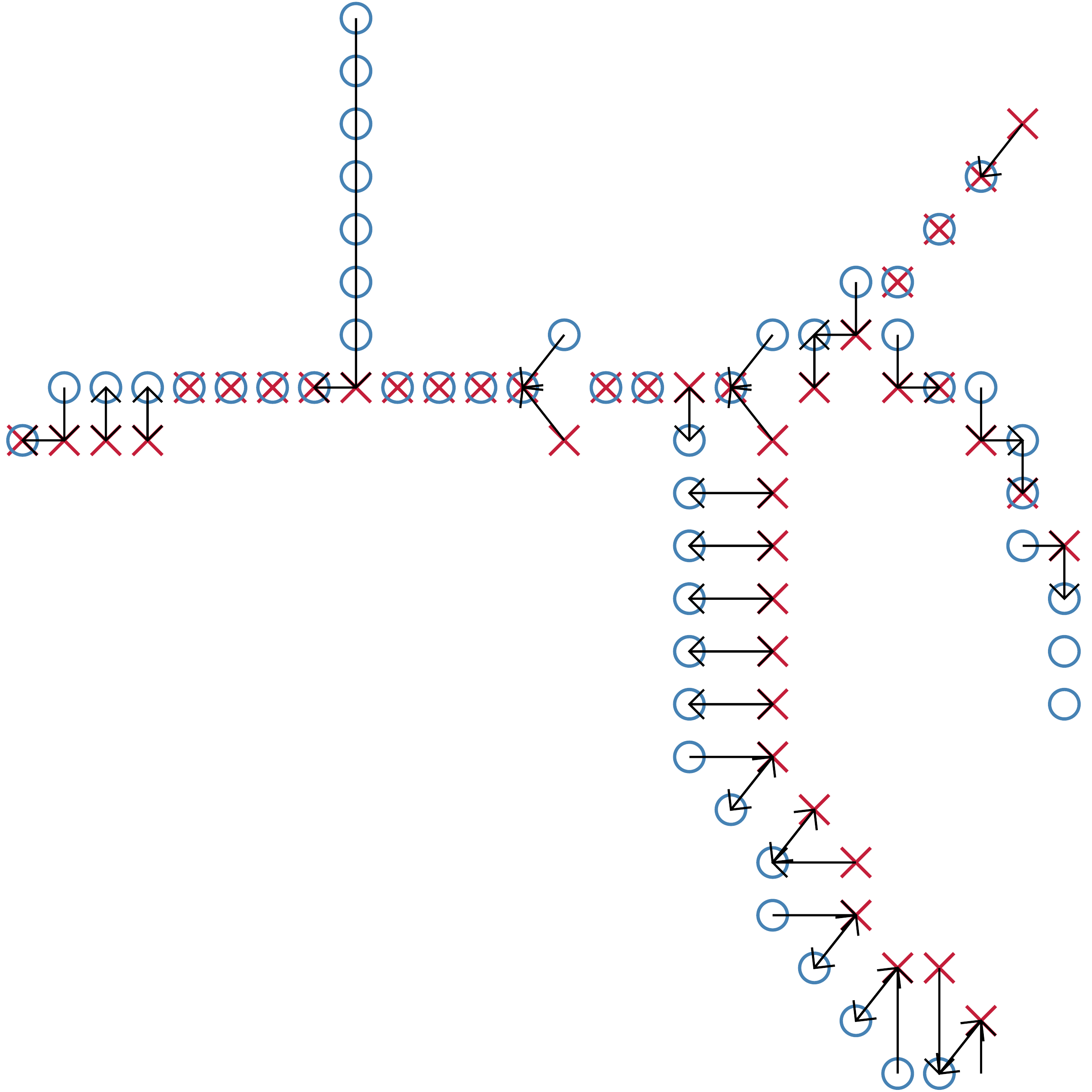}
  }
  \subfigure[$t=74$]{%
    \includegraphics[height=3.75cm]{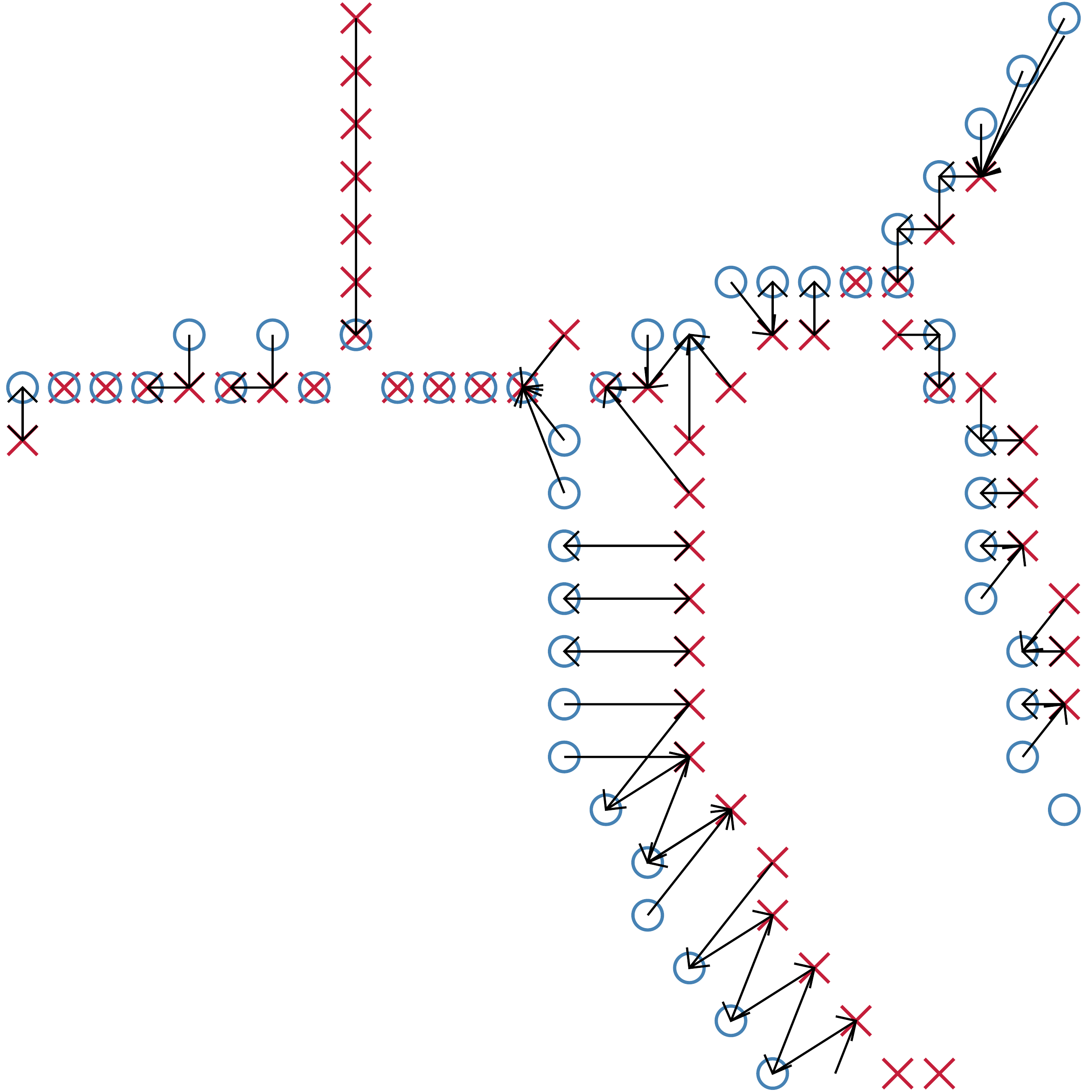}
  }
  \subfigureCaptionSkip
  \caption{%
    An excerpt demonstrating matches between two time steps. Some of the pixels of the current time
    step~(blue circles) overlap with pixels from the previous time step~(red crosses). We use arrows to indicate forward and backward matches.
  }
  \label{fig:Matches example}
\end{figure}

%
Since the skeleton changes over time, we need to characterize the creation time of each pixel, i.e.,
the time step~$t_i$ in which it initially appears.
Moreover, we want to permit that a pixel ``moves'' slightly between two consecutive time steps in order
to ensure that drifts of the skeleton can be compensated.
Our experiments indicate that it is possible to obtain consistent creation times for the pixels
based on their nearest neighbors, regardless of whether the simulation suffers from a coarse time
resolution or not. Given two time steps $t_i$, $t_{i+1}$, we assign every pixel $p
\in \pixels_i$ the pixel~$p' \in \pixels_{i+1}$ that satisfies
\begin{equation}
  p' := \arg\min_{q\in\pixels_{i+1}} \dist(p, q),
\end{equation}
where $\dist(\cdot)$ is the Euclidean distance. Likewise, we assign every pixel in $\pixels_{i+1}$
its nearest neighbor in $\pixels_{i}$, which represents a match from $\pixels_{i+1}$ to $\pixels_i$.
This yields a set of directed matches between~$\pixels_i$ and $\pixels_{i+1}$.  Each pixel is
guaranteed to occur at least once in the set.
We refer to matches from $\pixels_i$ to $\pixels_{i+1}$ as \emph{forward} matches, while we refer to
matches in the other direction as \emph{backward} matches. A match is \emph{unique} if the forward
and backward match connect the same pair of pixels.
Fig.~\ref{fig:Matches example} depicts matches for selected time steps and illustrates the movement
of pixels.

\subsubsection{Pixel Classification}
%
\begin{figure}[tbp]
  \centering%
  \subfigure[$t=68$]{%
    \includegraphics[width=0.32\linewidth]{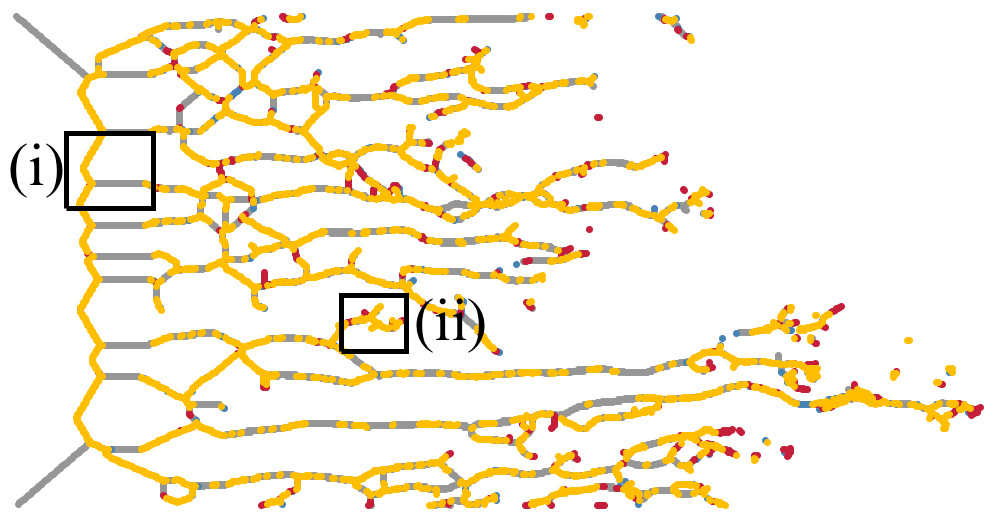}%
  }%
  \subfigure[$t=69$\label{sfig:t69 classification}]{%
    \includegraphics[width=0.32\linewidth]{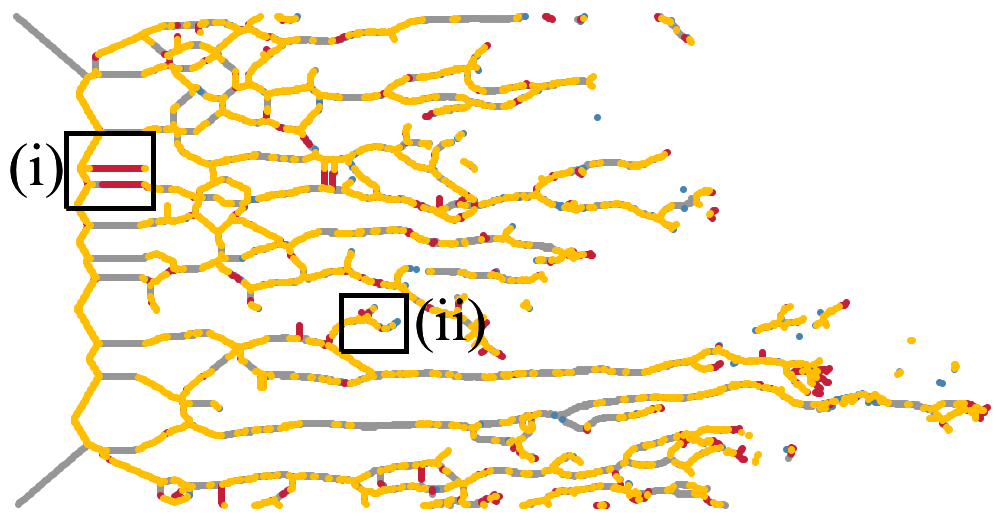}%
  }%
  \subfigure[$t=70$\label{sfig:t70 classification}]{%
    \includegraphics[width=0.32\linewidth]{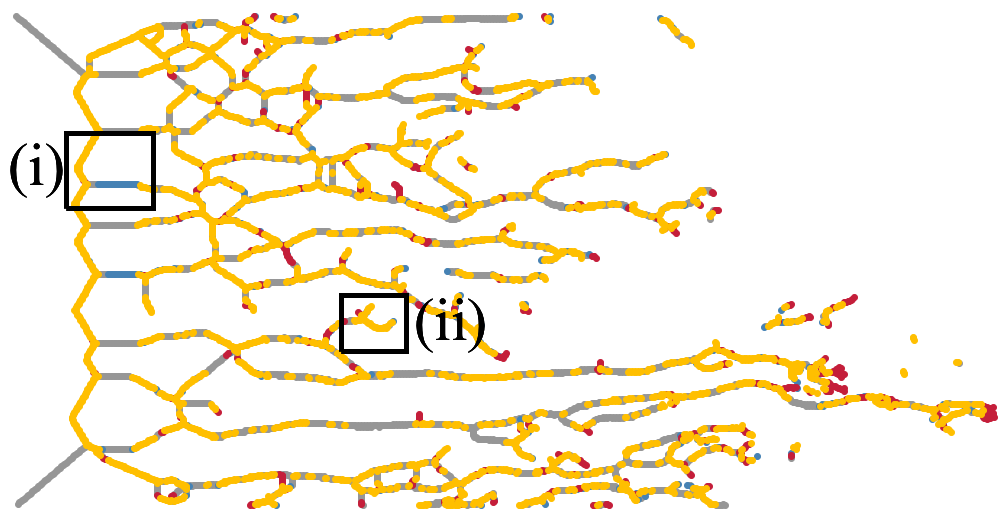}%
  }
  \subfigureCaptionSkip
  \caption{%
    Classification of all pixels into growth pixels~\inlinedot{created}, decay
    pixels~\inlinedot{destroyed}, known pixels~\inlinedot{persisting}, and irregular
    pixels~\inlinedot{irregular}.
    The abrupt appearance~\subref{sfig:t69 classification} or disappearance~\subref{sfig:t70
    classification} of segments is a challenge for skeleton extraction and tracking.
  }
  \label{fig:Pixel classification}
\end{figure}
%
We now classify each pixel in time step $t_{i+1}$ according to the forward matches between
$\pixels_i$ and $\pixels_{i+1}$, as well as the backward matches between $\pixels_i$ and
$\pixels_{i+1}$.
We call a pixel \emph{known} if their match is unique, i.e., there is exactly one forward and one
backward match that relate the same pixels with each other.
Known pixels are pixels that are already present in a previous time step with a unique counterpart
in time step~$t_{i+1}$.
Similarly, we refer to a pixel in $\pixels_{i+1}$ as a \emph{growth} pixel if there is a unique
match in $\pixels_i$ and at most one forward match from some other pixel in $\pixels_i$.
Growth pixels indicate that new structures have been created in time step $t_{i+1}$, or that
existing structures have been subject to a deformation.
The counterpart to a growth pixel is a \emph{decay} pixel in $\pixels_{i+1}$, which is defined by
a unique match in $\pixels_i$ and at most one backward match to the same pixel in $\pixels_i$ from
another pixel in $\pixels_{i+1}$.
Decay pixels indicate that a skeleton region has been lost in time step $t_{i+1}$.
We refer to all other pixels as \emph{irregular}.
In our experiments, irregular pixels, which are caused by small shifts between consecutive time
steps, comprise about 60\% of all pixels. As we subsequently demonstrate, we are able to assign
consistent creation times despite the prevalence of irregular pixels.
Fig.~\ref{fig:Pixel classification} depicts classified pixels for consecutive time steps. It
also demonstrates that skeletons may be temporally incoherent: pixels in region~(i) only
exist for a single time step, forming long but short-lived segments.  Pixels in region~(ii), by
contrast, form short but long-lived segments. We want to filter out segments in region~(i), while
keeping segments in region~(ii) intact. This requires knowledge about pixel creation times.

\subsubsection{Propagating Creation Times}
\label{sec:Propagation}
%
Initially, each pixel in $\pixels_0$ is assigned a creation time of $0$. Next, we classify the
pixels in each pair of consecutive time steps $t_i$ and $t_{i+1}$ as described above. For known pixels,
we re-use the creation time of $t_i$.
For growth pixels, we distinguish two different cases:
\begin{inparaenum}[(i)]
  \item If a growth pixel in time step $t_{i+1}$ is not the target of a forward match from time step
  $t_i$, we consider it to be a new pixel and hence assign it a creation time of $t_{i+1}$.
  \item Else, we re-use the creation time just as for known pixels.
\end{inparaenum}
This procedure ensures that we are conservative with assigning ``new'' creation times; it turns out
that a small number of growth pixels with increased creation times is sufficient for propagating
time information throughout the data.
For all other types of pixels, we assign them the minimum of all creation times of their respective
matches from $\pixels_i$, ignoring the direction of the matching. Again, this is a conservative
choice that reduces the impact of noise in the data.

\begin{figure}[tbp]
  \centering%
  \subfigure[$t=68$]{%
    \includegraphics[width=0.32\linewidth]{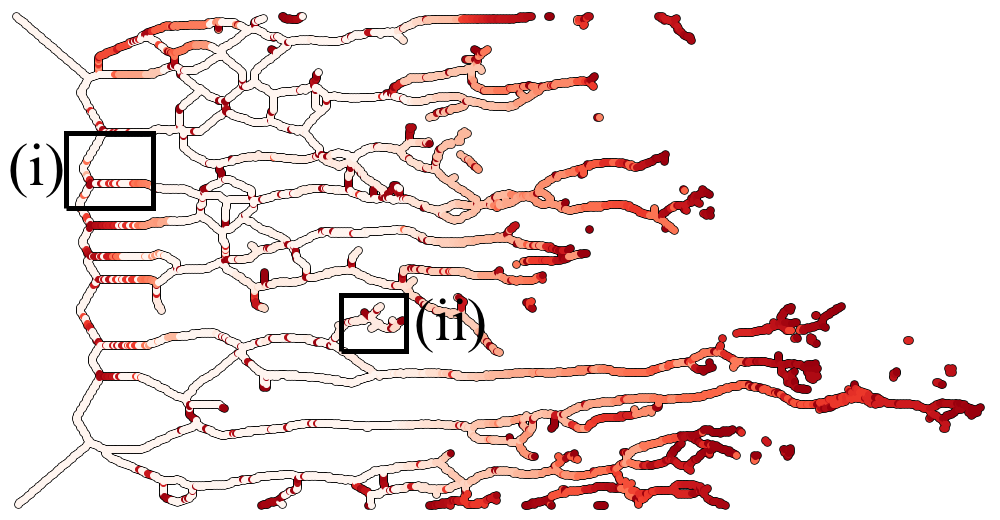}%
  }%
  \subfigure[$t=69$\label{sfig:t69 age}]{%
    \includegraphics[width=0.32\linewidth]{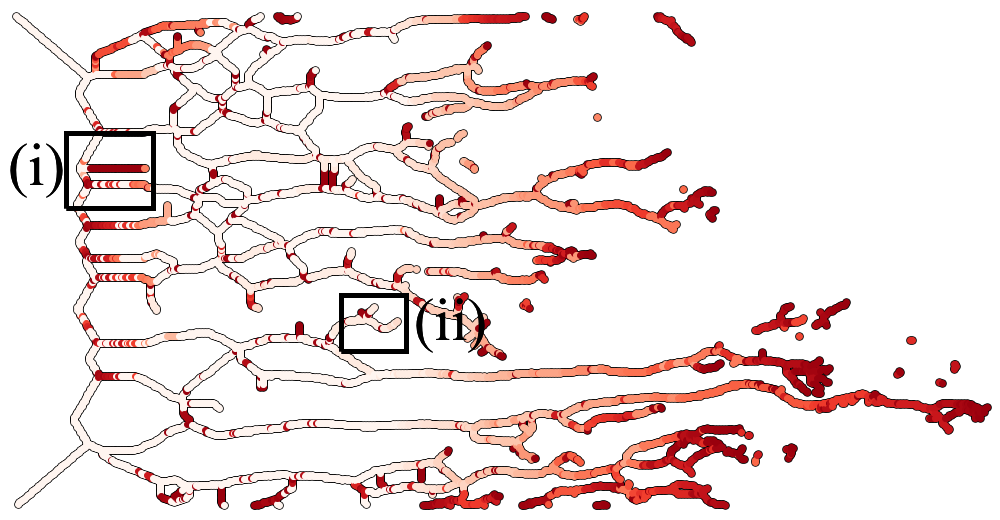}%
  }%
  \subfigure[$t=70$\label{sfig:t70 age}]{%
    \includegraphics[width=0.32\linewidth]{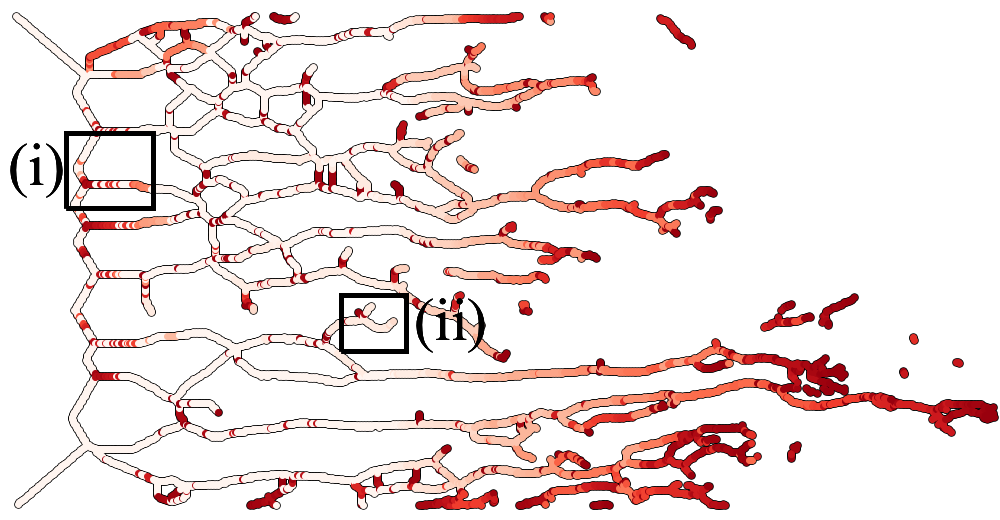}%
  }
  \subfigureCaptionSkip
  \caption{%
    Propagated age per pixel, using a white--red color map. The skeleton inconsistencies in region~(i) impede the temporal coherence of neighboring pixels.
  }
  \label{fig:Pixel age}
\end{figure}

Thus, every pixel in every time step has been assigned a creation time.  This time refers to the
first time step in which the pixel was unambiguously identified and appeared.
By propagating the creation time, we ensure that skeletons are allowed to exhibit some
movement between consecutive time steps.
Fig.~\ref{fig:Pixel age} depicts the creation times for several time steps. For temporally coherent
skeletons, recent creation times~(shown in red) should only appear at the end of new ``fingers''. We
can see that the brief appearance of segments causes inconsistencies. Ideally, the creation time of
pixels should vary continuously among a segment.

\begin{figure}[tbp]
  \centering
  \subfigure[$t=42$~(no coherence)\label{sfig:Creation times t42 no coherence}]{%
    \includegraphics[height=3cm]{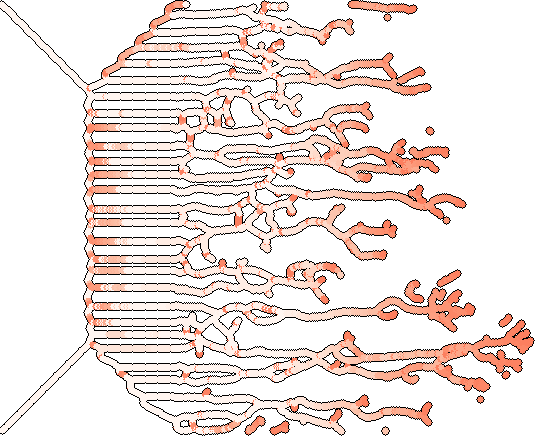}
  }
  \quad
  \subfigure[$t=84$~(no coherence)\label{sfig:Creation times t84 no coherence}]{%
    \includegraphics[height=3cm]{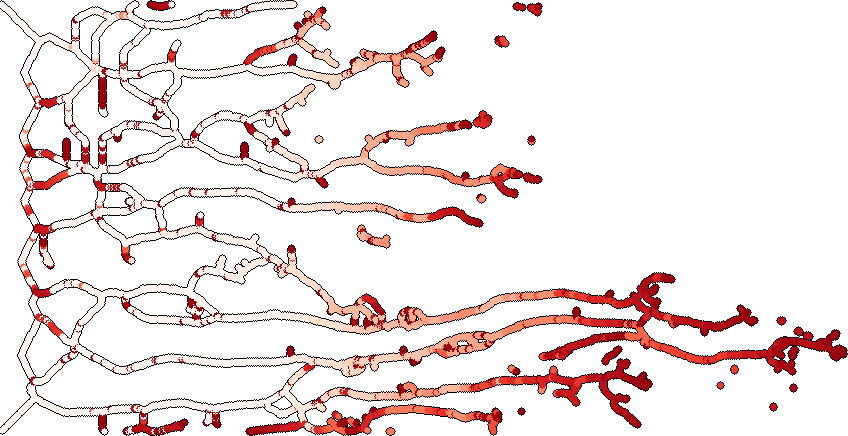}
  }
  \subfigure[$t=42$~(coherence)]{%
    \includegraphics[height=3cm]{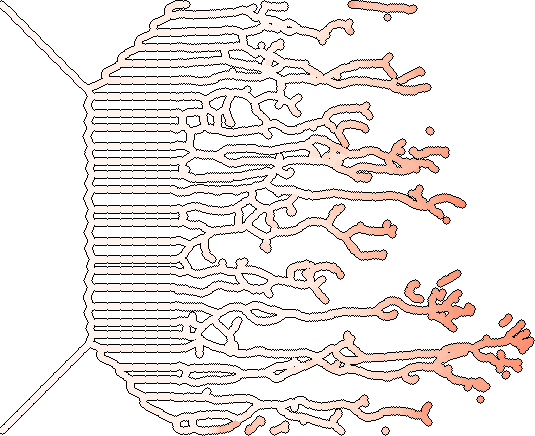}
  }
  \quad
  \subfigure[$t=84$~(coherence)]{%
    \includegraphics[height=3cm]{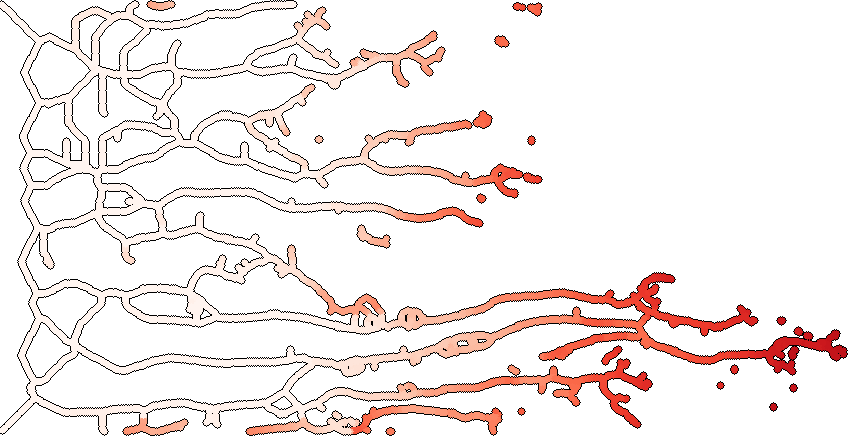}
  }
  \subfigureCaptionSkip
  \caption{%
    Pixel creation times at two selected time steps. Recent creation times are shown in shades of
    red. We can see that the ``front'' of the fingers is always recent, while the oldest structures have been created at the very beginning.
    This example also demonstrates how the temporal coherence of creation times can be improved.
  }
  \label{fig:Creation times example}
\end{figure}

\runinhead{Improving temporal coherence}
%
To improve temporal coherence, i.e, creation times of adjacent pixels, we observe that
inconsistencies are mainly caused by a small number of growth pixels along a segment. These are
a consequence of a ``drift'' in pixel positions over subsequent time steps, which our naive matching algorithm cannot compensate for.
A simple neighborhood-based strategy is capable of increasing coherence, though: for each growth
pixel, we evaluate the creation times in its $8$-neighborhood. If more than 50\% of the neighbors
have a different creation time than the current pixel, we replace its creation time by the
\emph{mode} of its neighbors' creation times. This strategy is reminiscent of \emph{mean shift
smoothing}~\cite{Comaniciu02}.
Fig.~\ref{fig:Creation times example} compares the original and improved creation times for two time
steps. Ideally, all segments should exhibit a gradient-like behavior, indicating that their
structures have been expanded continuously. We see that this is only true for the longest segments.
Erroneous creation times are an inevitable byproduct of instabilities in skeleton extraction, which
can be mitigated through persistence-based concepts. 

\subsection{Persistence Concepts}
\label{sec:Persistence concepts}
%
Persistence is a concept introduced by Edelsbrunner et
al.~\cite{Edelsbrunner10,Edelsbrunner02,Edelsbrunner06}.
It yields a measure of the range~(or scale) at which topological features occur in data and is
commonly employed to filter or simplify complex multivariate data sets~\cite{Rieck12}.
For skeletons, i.e., \emph{graphs}, the standard topological features are well known, comprising
connected components and cycles. While these features are useful in classifying complex
networks~\cite{Carstens13}, for example, they do not provide sufficient information about skeleton
evolution processes because they cannot capture the growth of segments.
Hence, instead of adopting this viewpoint, we derive several concepts that are inspired by the
notion of persistence. A crucial ingredient for this purpose is the availability of creation times
for every pixel in every time step.

\subsubsection{Branch Inconsistency}
%
Using the graph~$\graph_i$ for a time step $t_i$, we know which pixels are branch points, i.e.,
points where multiple segments meet.
Let $c_b$ be the creation time of such a branch point, and let $c_1$, $c_2$, \dots refer to the
creation times of the first adjacent point along each of the segments meeting at the branch point.
We define the \emph{branch inconsistency} for each branch--segment pair as~$|c_i - c_b|$, and we refer
to the diagram formed by the points $(c_b, c_i)$ as the branch persistence diagram.
The number of points in the branch inconsistency diagram indicates how many new branches are being
created in one time step.
Moreover, it can be used to prune away undesired segments in a skeleton: if the branch inconsistency
of a given segment is large, the segment is likely an artifact of the skeletonization
process---thinning algorithms often create segments that only exist for a single time step. Overall,
those segments thus have a late creation time.
In contrast to the persistence diagrams in topological data analysis, where closeness to the
diagonal indicates noise, here, points that are \emph{away} from the diagonal correspond to
erroneous segments in the data. Points \emph{below} the diagonal are the result of
inconsistent creation times for some segments---a branch cannot be created \emph{before} its
branch point.

\begin{figure}[tbp]
  \centering
  \subfigure[Branch inconsistency]{%
    \includegraphics[width=3cm]{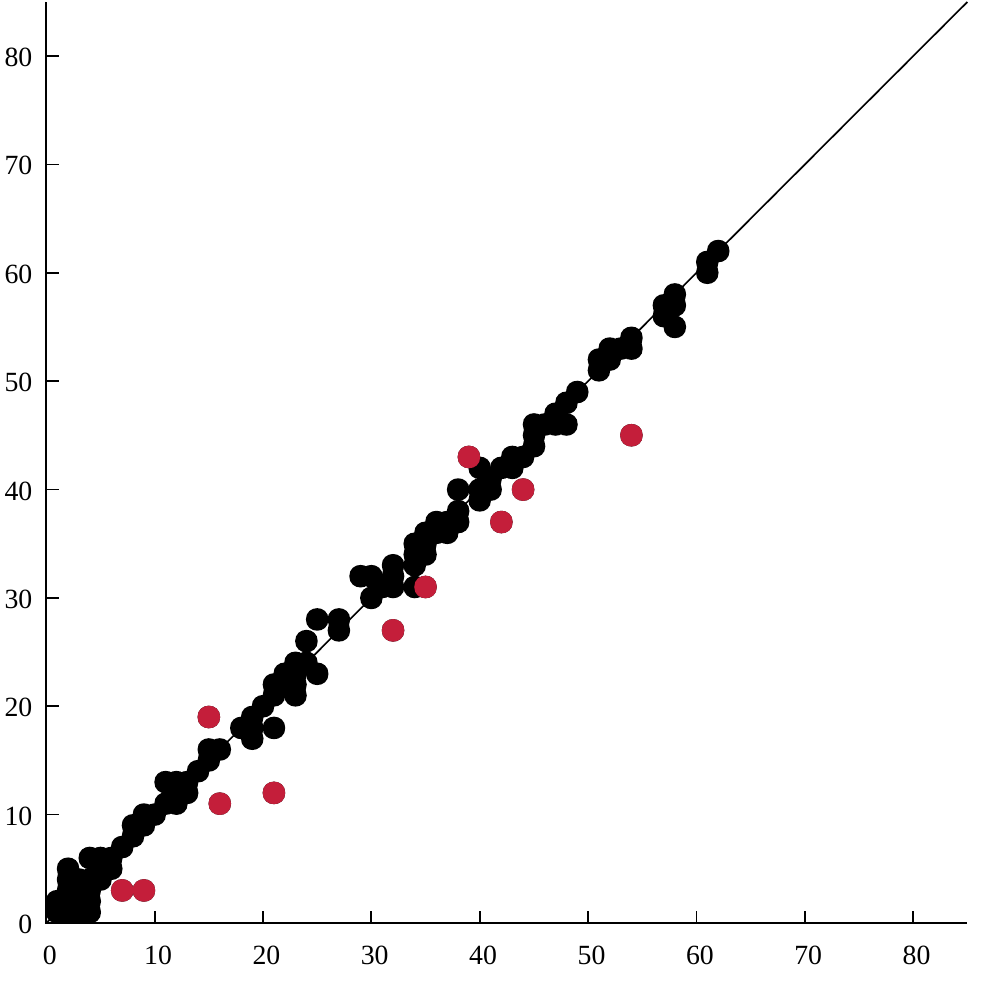}
  }%
  \subfigure[Skeleton]{%
    \raisebox{0.30\height}{%
      \includegraphics[width=0.33\linewidth]{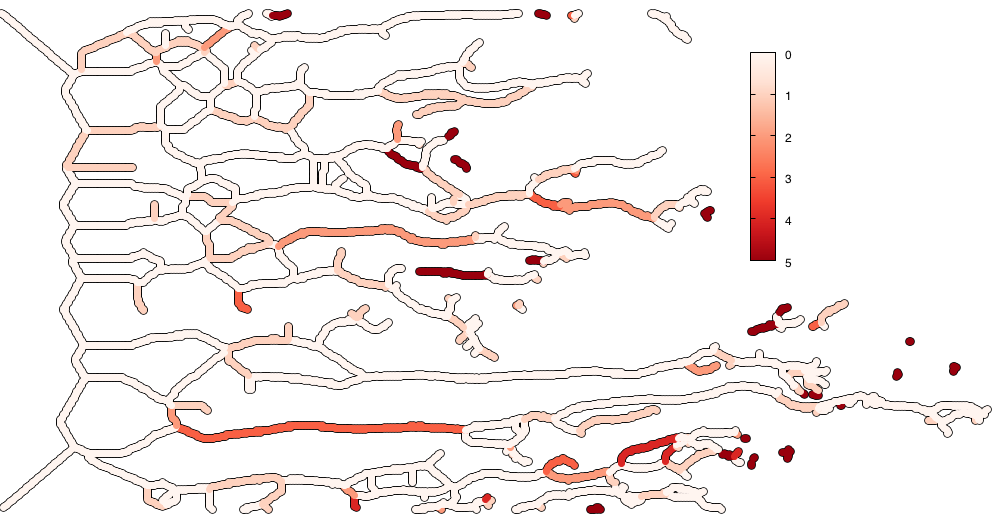}%
    }
  }%
  \subfigure[Skeleton, filtered\label{sfig:Branch persistence filtered}]{%
    \raisebox{0.30\height}{%
      \includegraphics[width=0.33\linewidth]{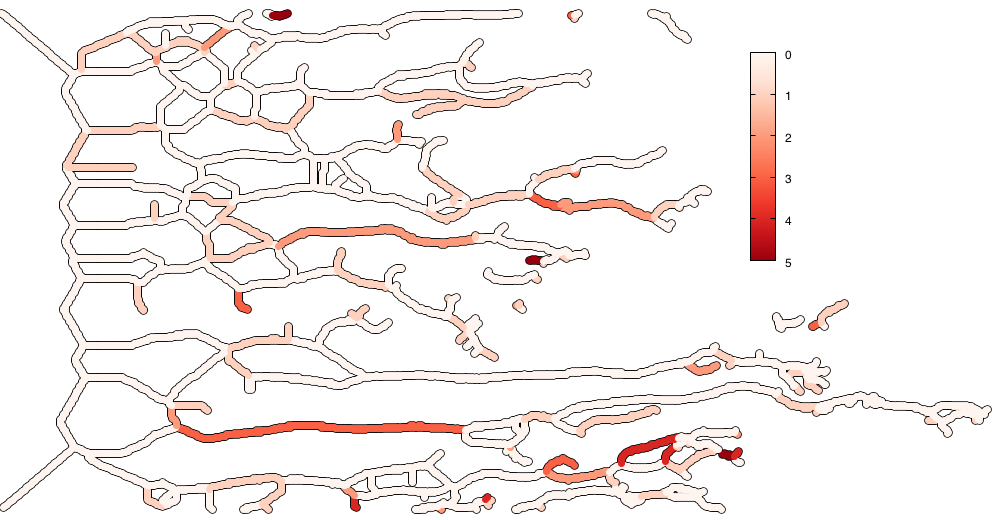}
    }
  }
  \subfigureCaptionSkip
  \caption{%
    Branch inconsistency diagram and branch inconsistency values on the skeleton for $t=69$.
    The diagram indicates that most branches are temporally coherent. Some of them are removed from
    the diagonal~(or below the diagonal), which may either indicate inconsistencies in skeleton
    tracking or cycles. Removing segments with a branch inconsistency $\geq 5$~(red dots in the
    diagram, dark red segments in the skeleton) can be used to filter the skeleton.
  }
  \label{fig:Branch persistence}
\end{figure}

Fig.~\ref{fig:Branch persistence} shows the branch inconsistency diagram and colored skeletons for
$t=69$. It also depicts how to filter segments with a large branch inconsistency, which already
decreases the number of noisy segments.
Please refer to the accompanying video for all branch inconsistency values.

\subsubsection{Age Persistence}
%
Analogously to branch inconsistency, we obtain an \emph{age persistence} diagram for each branch--segment pair when we use
the maximum creation time of points along each segment. Age persistence is capable of measuring
whether a segment is young or old with respect to its branch point.
Here, the ``persistence'' of each point is an indicator of how much the skeleton grows over multiple
time steps: if segments stagnate, their points remain at the same distance from the diagonal. If segments continue to grow, however, their points will move away from the diagonal.

\begin{figure}[tbp]
  \centering
  \subfigure[Age persistence]{%
    \includegraphics[width=3cm]{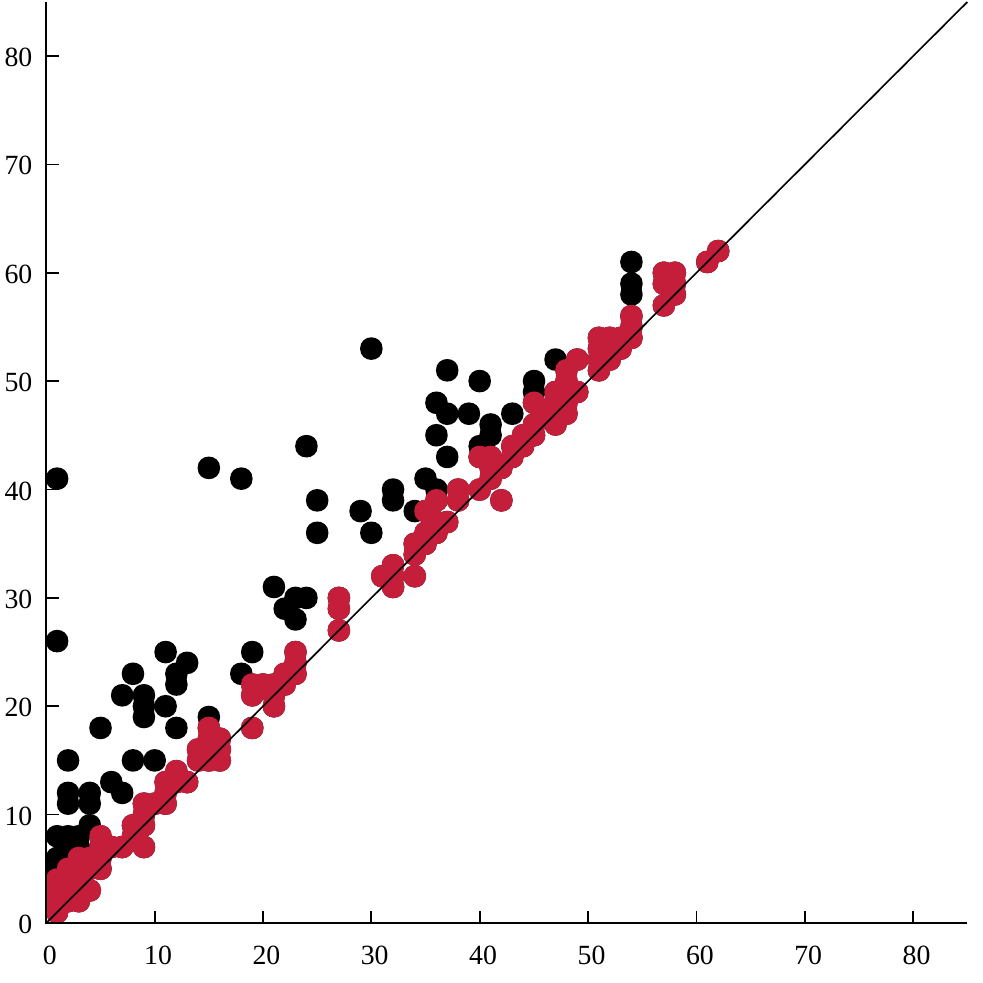}
  }%
  \subfigure[Skeleton]{%
    \raisebox{0.30\height}{%
      \includegraphics[width=0.33\linewidth]{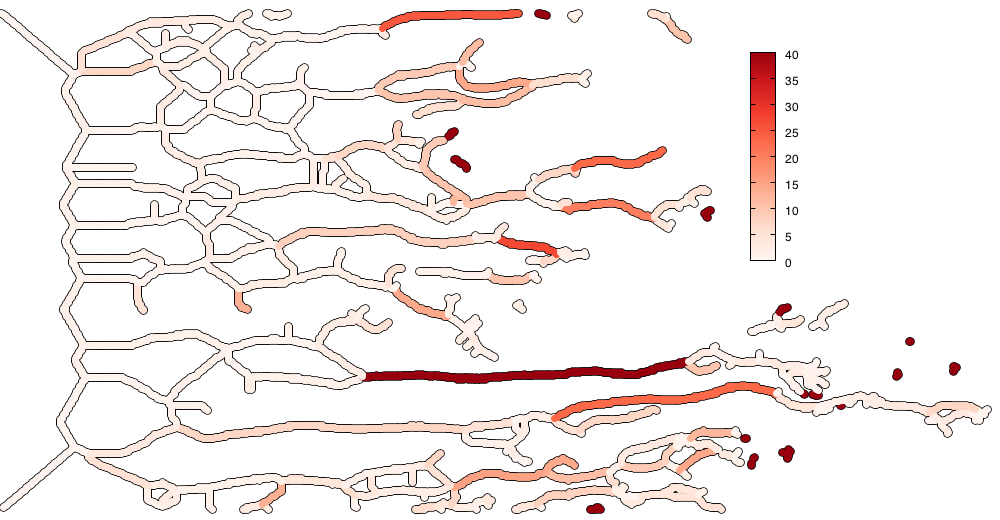}%
    }
  }%
  \subfigure[Skeleton, filtered]{%
    \raisebox{0.30\height}{%
      \includegraphics[width=0.33\linewidth]{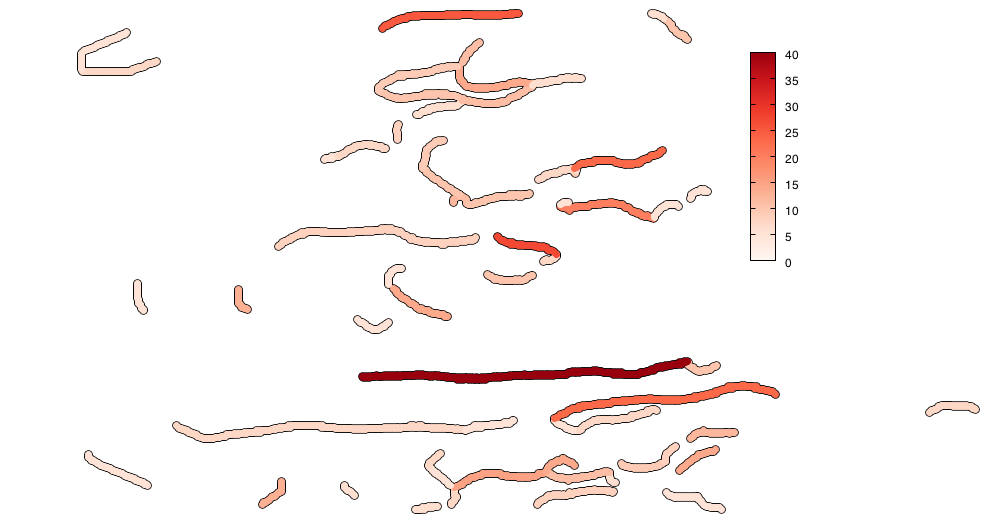}%
    }
  }
  \subfigureCaptionSkip
  \caption{%
    Age persistence diagram and age persistence values on the skeleton for $t=69$.
    Numerous segments towards the ``front'' of the fingers appear to be active here.
    Removing all segments whose age persistence is $\leq 5$~(red points in the diagram) leaves us
    with the most active segments.
  }
  \label{fig:Age persistence}
\end{figure}

\begin{figure}[btp]
  \centering
  \subfigure[\label{sfig:Branch_consistency}]{%
    \includegraphics[width=5cm]{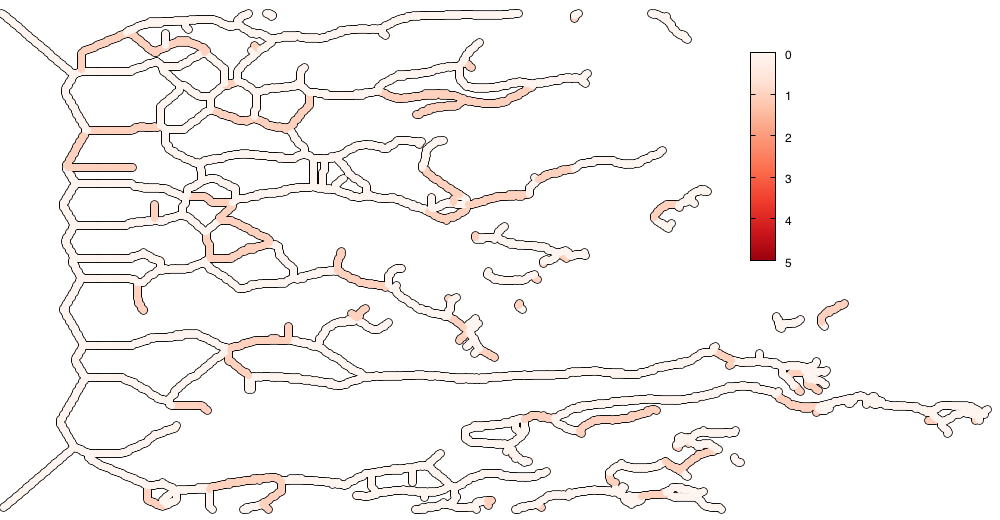}%
  }
  \subfigure[\label{sfig:Branch_consistency_and_age_persistence}]{%
    \includegraphics[width=5cm]{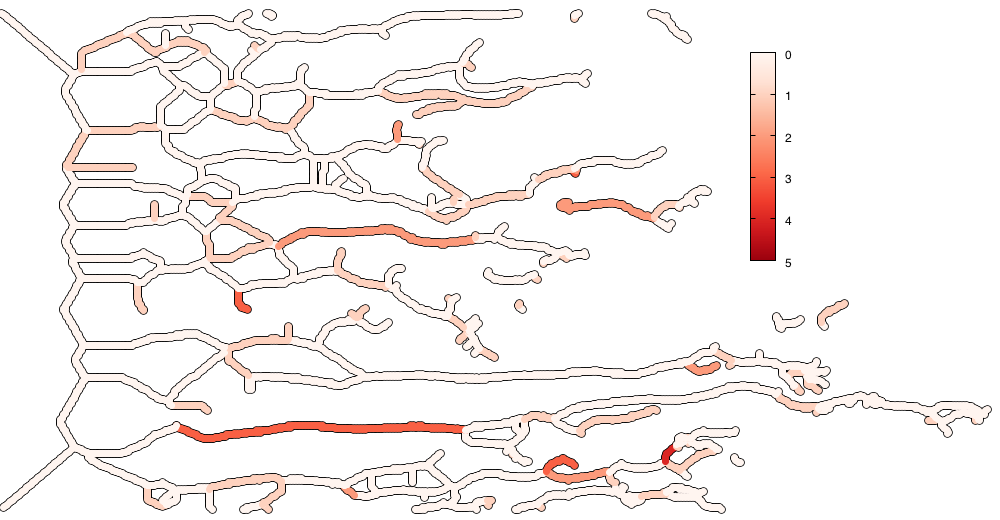}%
  }
  \subfigureCaptionSkip
  \caption{%
    \protect\subref{sfig:Branch_consistency} Filtering segments using branch inconsistency may
    destroy longer segments.
    \protect\subref{sfig:Branch_consistency_and_age_persistence} If we combine both branch
    inconsistency and age persistence, keeping only those segments whose age persistence is high or
    whose branch inconsistency is low, we can improve the filter results by removing noisy segments
    while keeping more cycles intact.
  }
  \label{fig:Combined filtering}
\end{figure}

Fig.~\ref{fig:Age persistence} shows the age persistence diagram and the age persistence values on
the skeleton for $t=69$. The filtered skeleton only contains the most active segments, which
facilitates tracking.
We can combine branch inconsistency and age persistence to remove fewer segments than in
Fig.~\ref{sfig:Branch persistence filtered}. For example, we could remove segments that
correspond to points below the diagonal of the branch inconsistency diagram and keep those for which
\emph{both} branch inconsistency and age persistence are high. These segments commonly correspond to
cycles that were formed during the evolution of the skeleton. An isolated analysis of branch
inconsistency is unable to detect them. Fig.~\ref{fig:Combined filtering} depicts the results of
such a combined filtering operation.

\subsubsection{Growth Persistence}
\begin{figure}[tbp]
  \centering
  \subfigure[$t=21$]{\includegraphics[height=2.5cm]{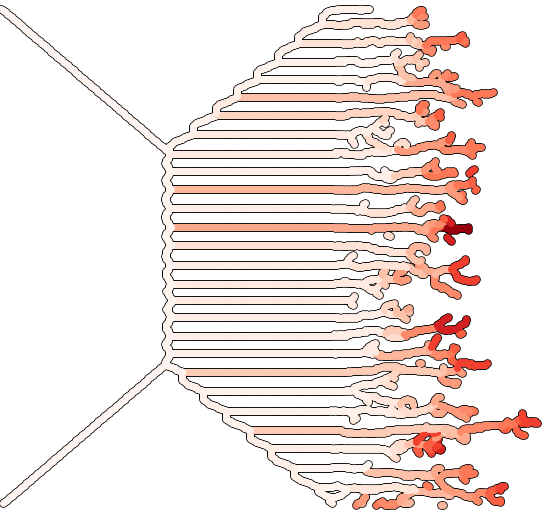}}
  \subfigure[$t=42$]{\includegraphics[height=2.5cm]{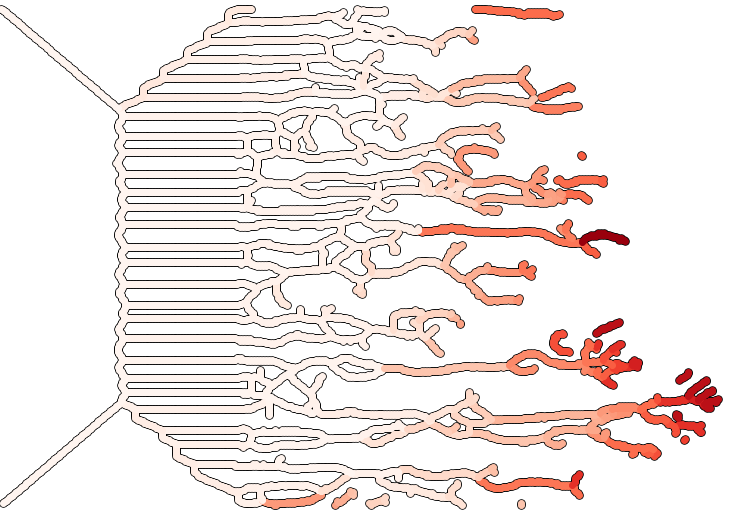}}
  \subfigure[$t=84$]{\includegraphics[height=2.5cm]{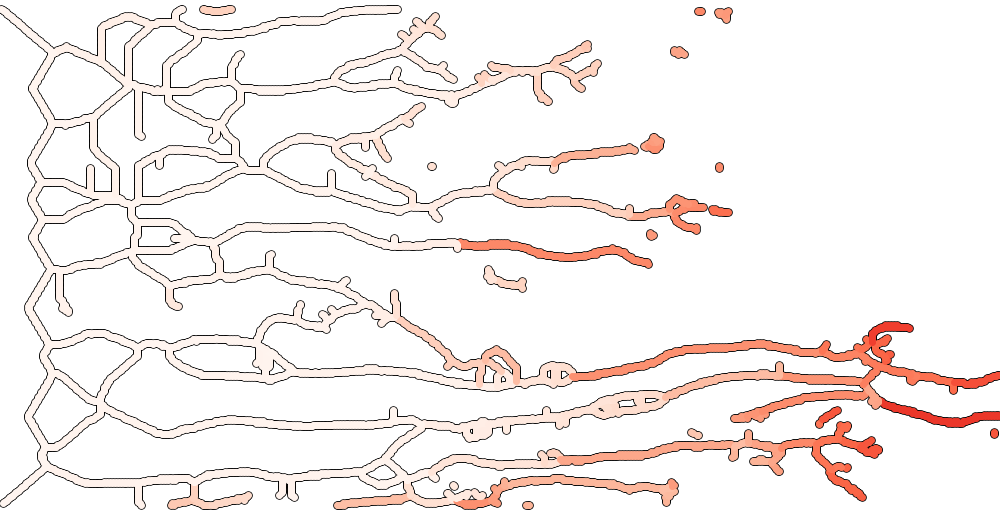}}
  \subfigureCaptionSkip
  \caption{%
    Growth persistence values. Red segments are highly active in the evolution of the skeleton. In
    this example, red segments are mostly those that are at the tips of individual ``fingers''.
  }
  \label{fig:Growth persistence}
\end{figure}

%
We define the \emph{growth persistence} of a segment in $\graph_i$ as the difference between the
maximum creation time~$t_{\max}$ of its pixels and the current time step~$t_i$. Intuitively, this
can be thought of performing ``time filtration'' of a simplicial complex, in which simplices may be
created \emph{and} destroyed~(notice that such a description would require zigzag persistence for
general simplicial complexes).
A small value in this quantity indicates that the segment is still growing, while larger values
refer to segments that stagnate.
Growth persistence is useful to highlight segments that are relevant for tracking in viscous
fingering processes.
In contrast to the previously-defined persistence concepts, growth persistence is only defined per
segment and does not afford a description in terms of a persistence diagram.
Fig.~\ref{fig:Growth persistence} depicts the growth persistence of several time steps. Red segments
are growing fast or have undergone recent changes, such as the creation of cycles. A low branch
persistence in segments, coupled with a low growth persistence corresponds to features that are
``active'' during skeleton evolution.
Please refer to the accompanying video for the evolution of growth persistence.

\subsection{Activity Indicators}
\label{sec:Activity indicators}

In order to capture the dynamics of skeleton evolution, we require a set of activity indicators.
They are based on the previously-defined concepts and can be used to quickly summarize a time series
of evolving skeletons.

\subsubsection{Total Persistence}

There are already various summary statistics for persistence diagrams. The \emph{2\nd-order total
persistence}~$\pers(\diagram)$~\cite{Cohen-Steiner10} of a persistence diagram~$\diagram$ is defined
as
\begin{equation}
  \pers(\diagram)_2 := \left(\sum_{(c,d)\in\diagram} \pers^2(c,d)\right)^{\frac{1}{2}},
\end{equation}
i.e., the sum of powers of the individual persistence values~(i.e.\ coordinate differences) of the
diagram. Total persistence was already successfully used to assess topological activity in
multivariate clustering algorithms~\cite{Rieck16a}.

Here, the interpretation of total persistence depends on the diagram for which we compute it. Recall
that in a branch inconsistency diagram, points of high ``persistence'' indicate inconsistencies in
branching behavior.
Total persistence thus helps detect anomalies in the data; see Fig.~\ref{fig:Branch persistence
comparison} for a comparison of total branch
persistence in different data sets.
For age persistence, by contrast, high persistence values show that a skeleton segment is still
actively changing. The total age persistence hence characterizes the dynamics of the data,
e.g., whether many or few segments are active at each time step. Fig.~\ref{fig:Age persistence
comparison} depicts a comparison of total age persistence in different data sets.

\subsubsection{Vivacity}
\label{sec:Vivacity}
%
We also want to measure the ``vivacity'' of a viscous fingering process. To this end, we employ the
growth persistence values. Given a growth threshold $t_\mathrm{G}$, we count all growth pixels with
$\persG \leq t_\mathrm{G}$ and divide them by the total number of pixels in the given time step.
This yields a measure of how much ``mass'' is being created at every time step of the process.
Similarly, we can calculate vivacity based on \emph{segments} in the data. However, we found that
this does not have a significant effect on the results, so we refrain from showing the resulting
curves.
Fig.~\ref{fig:Vivacity example} depicts vivacity curves for different data sets with $t_\mathrm{G} = 10$.

\section{Analysis}
\label{sec:Analysis}

\begin{figure}[tbp]
  \centering
  \subfigure[Measured data set, $t=33$\label{sfig:Measured data}]{%
    \includegraphics[angle=90,height=3.25cm]{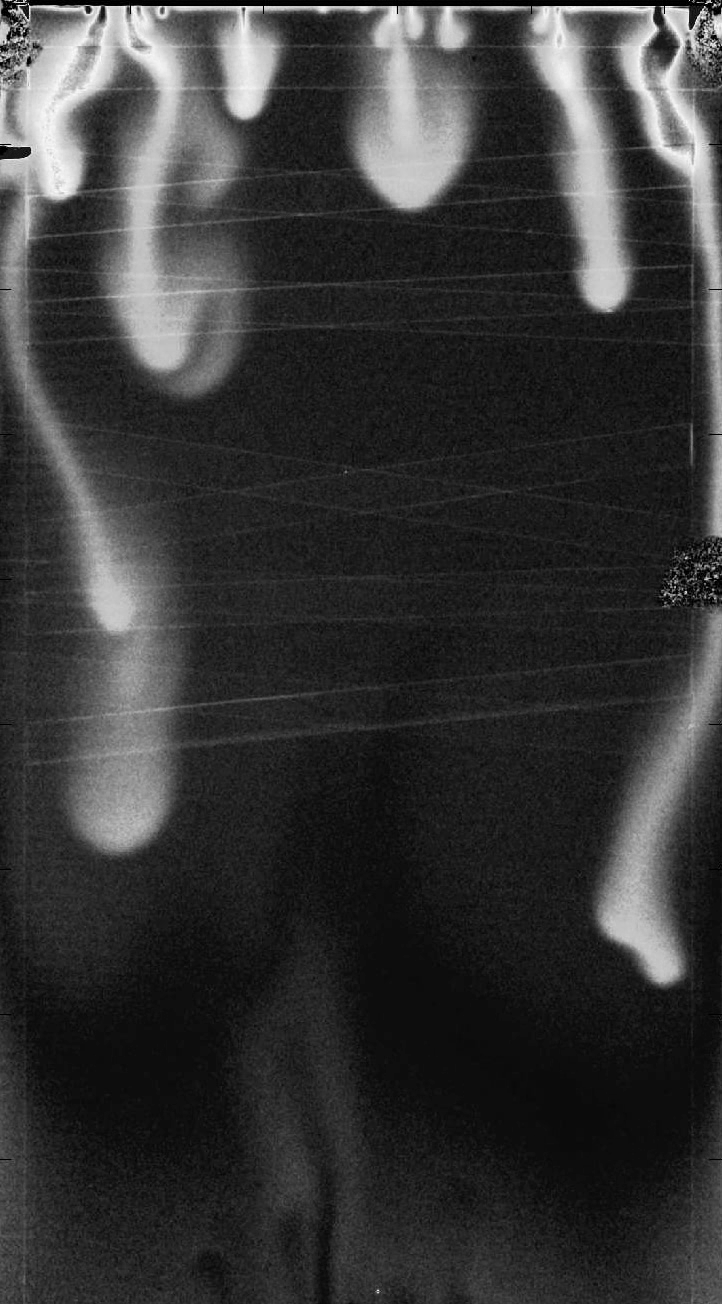}
  }
  \subfigure[Simulation data set, $t=26$\label{sfig:Simulation data}]{%
    \includegraphics[height=3.25cm]{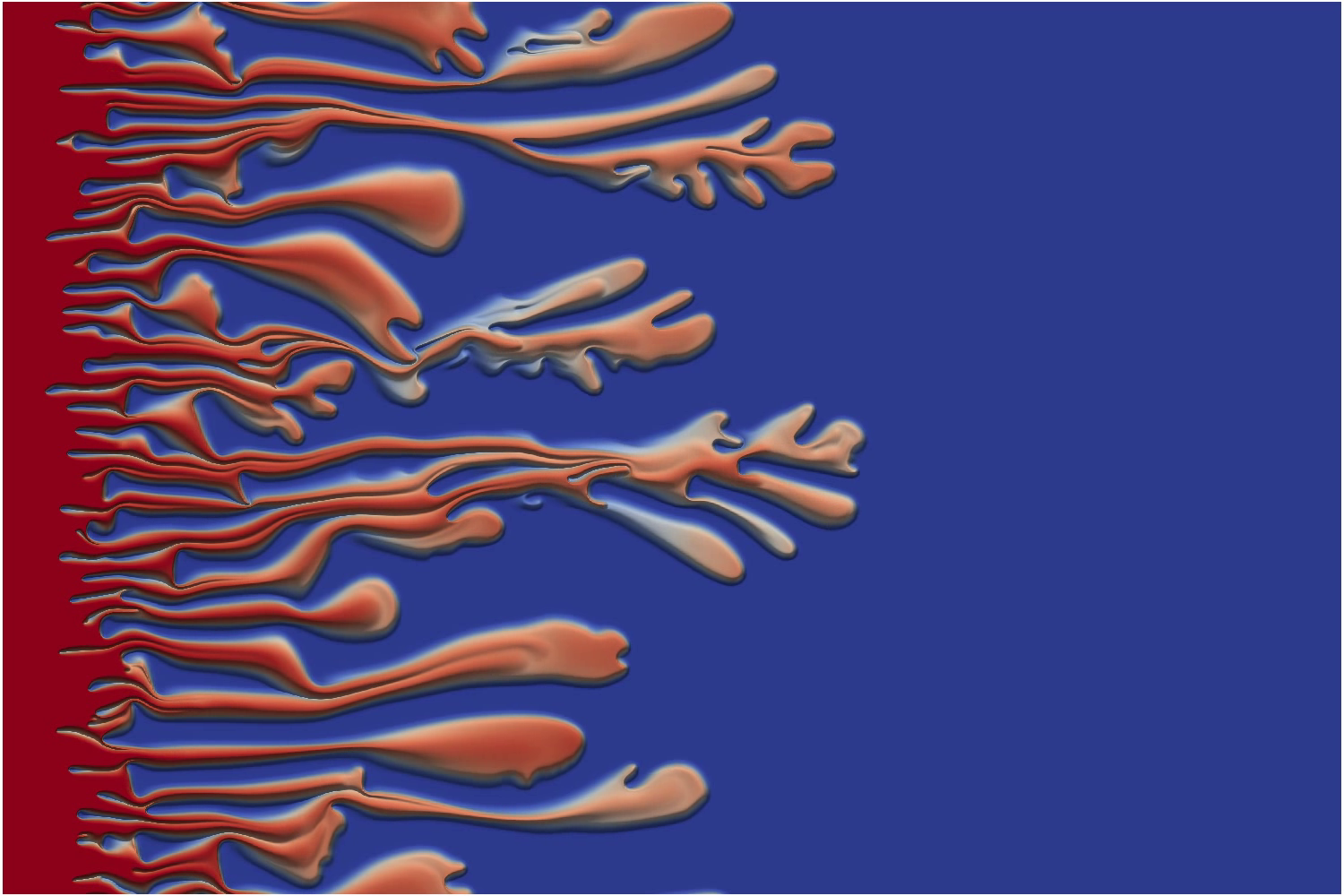}
  }
  \subfigureCaptionSkip
  \caption{%
    Selected still images from the two remaining data sets. The measured data
    in~\subref{sfig:Measured data} exhibits artifacts~(parallel lines) that are caused by  the
    experimental setup. The simulation data~\subref{sfig:Simulation data}, by contrast, does not contain any noise.
  }
  \label{fig:Data sets}
\end{figure}

Having defined a variety of persistence-based concepts, we now briefly discuss their utility in
analyzing time-varying skeleton evolution.
In the following, we analyze three different data sets:
\begin{inparaenum}[(i)]
  \item the \emph{example} data set that we used to illustrate all concepts,
  \item a \emph{measured} data set, corresponding to a slowly-evolving viscous fingering process,
  \item and a \emph{simulation} of the example data set.
\end{inparaenum}
Fig.~\ref{fig:Data sets} depicts individual frames of the latter two data sets.
The measured data set is characterized by a viscous fingering process whose fingers evolve rather
slowly over time. Moreover, this experiment, which was performed over several days, does not exhibit
many fingers. The simulation data, by contrast, aims to reproduce the dynamics found in the example
data set; hence, it contains numerous fast-growing fingers.
Please refer to the accompanying videos for more details.

\begin{figure}[tbp]
  \centering
  \subfigure[Example\label{sfig:Branch persistence example}]{%
    \iffinal
      \includegraphics{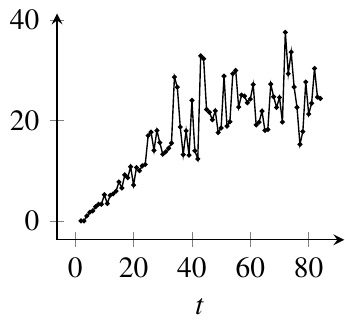}
    \else
      \begin{tikzpicture}
        \begin{axis}[%
          axis x line = bottom,
          axis y line = left,
          width       = 4.5cm,
          xlabel={$t$},
          enlargelimits = true,
        ]
          \addplot[sharp plot, mark=*, mark size=0.5pt] gnuplot[raw gnuplot] {%
            set datafile separator ',';
            set key autotitle columnheader;
            plot "Data/Skeletons-2/TXT/Matches-symmetric-branch-persistence-new/Branch_persistence_statistics.txt" using ($0+2):(sqrt($3));
          };
        \end{axis}
      \end{tikzpicture}
    \fi
  }
  \subfigure[Measured\label{sfig:Branch persistence 21_14}]{%
    \iffinal
      \includegraphics{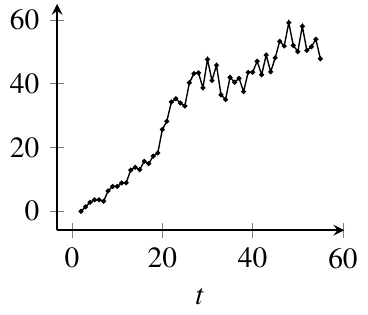}
    \else
      \begin{tikzpicture}
        \begin{axis}[%
          axis x line = bottom,
          axis y line = left,
          width       = 4.5cm,
          xlabel={$t$},
          enlargelimits = true,
        ]
          \addplot[sharp plot, mark=*, mark size=0.5pt] gnuplot[raw gnuplot] {%
            set datafile separator ',';
            set key autotitle columnheader;
            plot "Data/21_14/Branch-persistence/Branch_persistence_statistics.txt" using ($0+2):(sqrt($3));
          };
        \end{axis}
      \end{tikzpicture}
    \fi
  }
  \subfigure[Simulation\label{sfig:Branch persistence simulation}]{%
    \iffinal
      \includegraphics{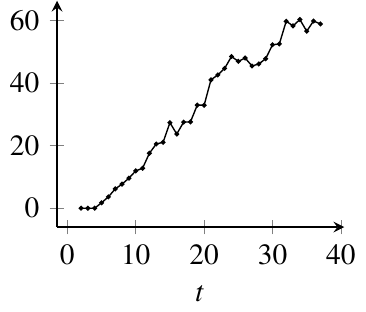}
    \else
      \begin{tikzpicture}
        \begin{axis}[%
          axis x line = bottom,
          axis y line = left,
          width       = 4.5cm,
          xlabel={$t$},
          enlargelimits = true,
        ]
          \addplot[sharp plot, mark=*, mark size=0.5pt] gnuplot[raw gnuplot] {%
            set datafile separator ',';
            set key autotitle columnheader;
            plot "Data/Simulation/Branch-persistence/Branch_persistence_statistics.txt" using ($0+2):(sqrt($3));
          };
        \end{axis}
      \end{tikzpicture}
    \fi
  }
  \subfigureCaptionSkip
  \caption{%
    A comparison of total persistence of the branch inconsistency diagram for three different data
    sets. The first data set~\subref{sfig:Branch persistence example} exhibits more
    anomalies; these are indicated by ``jumps'' in the total persistence curve.
  }
  \label{fig:Branch persistence comparison}
\end{figure}
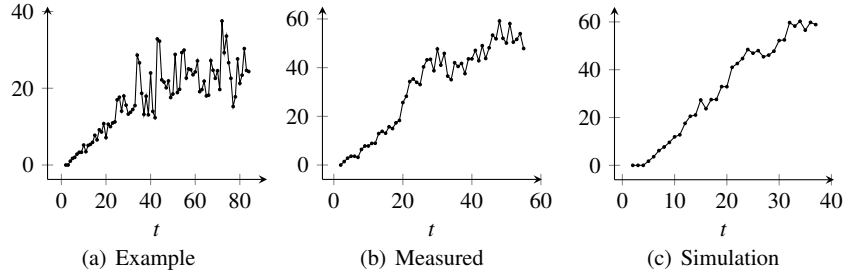

\begin{figure}[tbp]
  \centering
  \subfigure[$t=43$]{%
    \includegraphics[height=3.50cm]{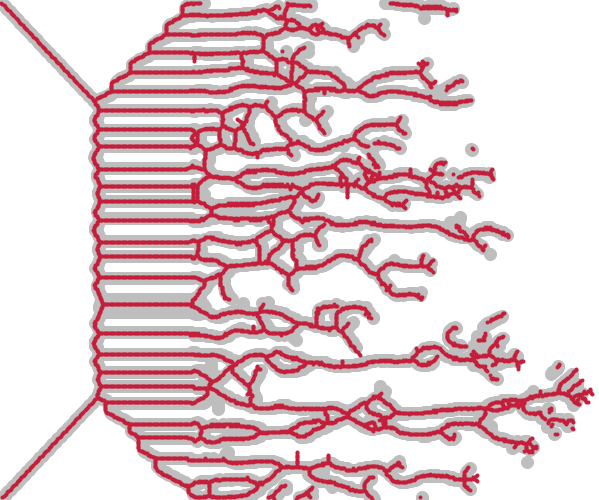}
  }
  \subfigure[$t=72$]{%
    \includegraphics[height=3.50cm]{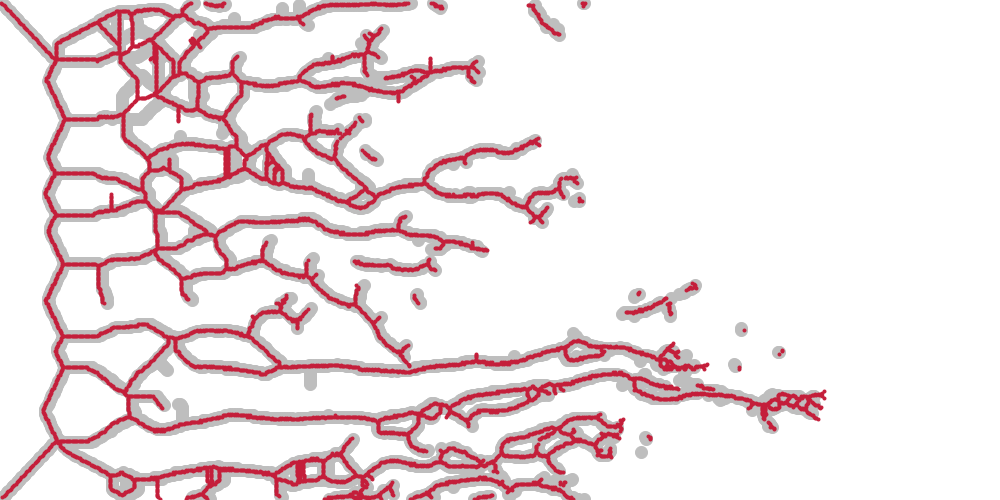}
  }
  \subfigureCaptionSkip
  \caption{%
    Comparing the previous~(gray) and the current time step~(red) based on total persistence of the
    branch inconsistency diagram helps uncover problems with skeleton extraction.
  }
  \label{fig:Anomaly detection example}
\end{figure}

\runinhead{Anomaly detection}
%
To detect anomalies in skeleton extraction and tracking, we calculate the total persistence of the
branch inconsistency diagram. Fig.~\ref{fig:Branch persistence comparison} compares the values for
all data sets.  We observe that the example data set, Fig.~\ref{sfig:Branch persistence example},
exhibits many ``jumps'' in branch inconsistency. These are time steps at which the skeleton~(briefly)
becomes inconsistent, e.g., because a large number of segments disappears, or many small cycles are
created.
At $t=43$ and $t=72$~(both local maxima in the diagram), for example, we observe changes in the
number of cycles as well as the appearance of numerous segments of various lengths, which makes it
harder to assign consistent creation times according to Sec.~\ref{sec:Extraction and propagation}.
Fig.~\ref{fig:Anomaly detection example} depicts the changes in skeleton topology at these time steps.
The other two data sets contain fewer anomalies. For the measured data, this is caused by lower
propagation velocities and fewer ``fingers'' in the data. For the simulation data, this is due to
a better separation of individual fingers, caused by the synthetic origin of the data.

\begin{figure}[tbp]
  \centering
  \subfigure[Example\label{sfig:Age persistence example}]{%
    \iffinal
      \includegraphics{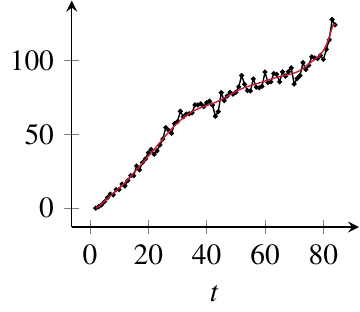}
    \else
      \begin{tikzpicture}
        \begin{axis}[%
          axis x line   = bottom,
          axis y line   = left,
          width         = 4.5cm,
          ytick         = {0,50,100},
          xlabel        = {$t$},
          enlargelimits = true,
          clip marker paths = true,
        ]

          \addplot[sharp plot, mark=*, mark size=0.5pt] gnuplot[raw gnuplot] {%
            set datafile separator ',';
            set key autotitle columnheader;
            plot "Data/Skeletons-2/TXT/Matches-symmetric-branch-persistence-new/Age_persistence_statistics.txt" using ($0+2):(sqrt($3));
          };

          \addplot[cardinal, sharp plot, mark=none] gnuplot[raw gnuplot] {%
            set datafile separator ',';
            set key autotitle columnheader;
            set samples 1000;
            plot "Data/Skeletons-2/TXT/Matches-symmetric-branch-persistence-new/Age_persistence_statistics.txt" using ($0+2):(sqrt($3)) smooth bezier;
          };

          %
          %
        \end{axis}
      \end{tikzpicture}
    \fi
  }
  \subfigure[Measured\label{sfig:Age persistence 21_14}]{%
    \iffinal
      \includegraphics{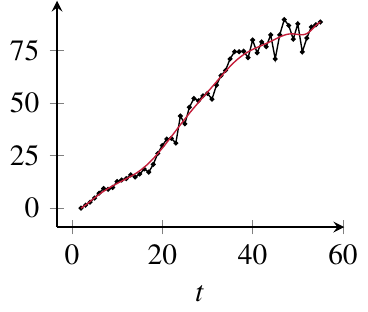}
    \else
      \begin{tikzpicture}
        \begin{axis}[%
          axis x line   = bottom,
          axis y line   = left,
          width         = 4.5cm,
          ytick         = {0,25,50,75,100},
          xlabel        = {$t$},
          enlargelimits = true,
          clip marker paths = true,
        ]
          \addplot[sharp plot, mark=*, mark size=0.5pt] gnuplot[raw gnuplot] {%
            set datafile separator ',';
            set key autotitle columnheader;
            plot "Data/21_14/Branch-persistence/Age_persistence_statistics.txt" using ($0+2):(sqrt($3));
          };
          \addplot[cardinal, sharp plot, mark=none] gnuplot[raw gnuplot] {%
            set datafile separator ',';
            set key autotitle columnheader;
            set samples 1000;
            plot "Data/21_14/Branch-persistence/Age_persistence_statistics.txt" using ($0+2):(sqrt($3)) smooth bezier;
          };
          %
          %
        \end{axis}
      \end{tikzpicture}
    \fi
  }
  \subfigure[Simulation\label{sfig:Age persistence simulation}]{%
    \iffinal
      \includegraphics{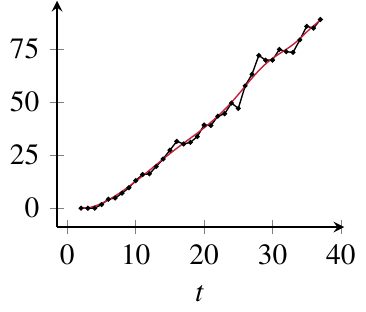}
    \else
      \begin{tikzpicture}
        \begin{axis}[%
          axis x line   = bottom,
          axis y line   = left,
          width         = 4.5cm,
          ytick         = {0,25,50,75,100},
          xlabel        = {$t$},
          enlargelimits = true,
        ]
          \addplot[sharp plot, mark=*, mark size=0.5pt] gnuplot[raw gnuplot] {%
            set datafile separator ',';
            set key autotitle columnheader;
            plot "Data/Simulation/Branch-persistence/Age_persistence_statistics.txt" using ($0+2):(sqrt($3));
          };
          \addplot[cardinal, sharp plot, mark=none] gnuplot[raw gnuplot] {%
            set datafile separator ',';
            set key autotitle columnheader;
            set samples 1000;
            plot "Data/Simulation/Branch-persistence/Age_persistence_statistics.txt" using ($0+2):(sqrt($3)) smooth bezier;
          };
        \end{axis}
      \end{tikzpicture}
    \fi
  }
  \subfigureCaptionSkip
  \caption{%
    A comparison of total age persistence for the three different data sets, along with a smooth
    estimate for showing trends.
  }
  \label{fig:Age persistence comparison}
\end{figure}
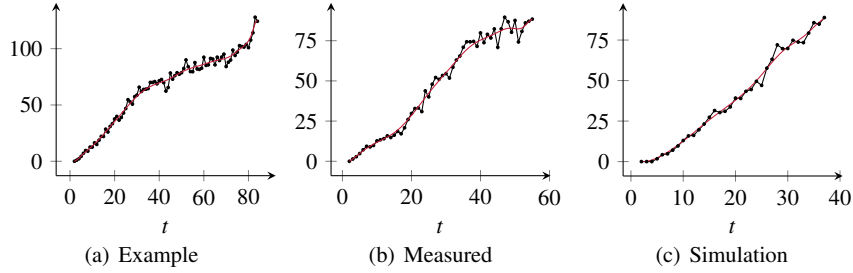

\runinhead{Active branches}
%
We use total age persistence to assess the rate at which existing branches move. Fig.~\ref{fig:Age
persistence comparison} compares the data sets, showing both the original total age persistence
values as well as a smooth estimate, obtained by fitting B{\'e}zier curves~\cite{Farin93} to the sample points.
In Fig.~\ref{sfig:Age persistence simulation}, the simulated origin of the data is evident: while
the other data sets exhibit changes in the growth rate of total age persistence, the simulation data
clearly exhibits almost constant growth.
Moreover, we observe that the measured data in Fig.~\ref{sfig:Age persistence 21_14} has a period of
constant growth for $t \in [20,40]$, while the example data displays a slightly diminished growth
rate for $t \in [25,65]$, only to pick up at the end. Age persistence may thus be used to compare
the characteristics of different skeleton evolution processes.

\begin{figure}[tbp]
  \centering
  \subfigure[Example\label{sfig:Vivacity pixel example}]{
    \iffinal
      \includegraphics{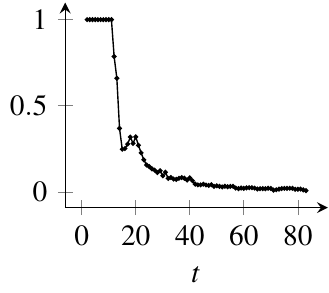}
    \else
      \begin{tikzpicture}
        \begin{axis}[%
          axis x line   = bottom,
          axis y line   = left,
          width         = 4.25cm,
          xlabel        = {$t$},
          enlargelimits = true,
        ]
          \addplot[sharp plot, mark=*, mark size=0.5pt] table[header=false] {Data/Vivacity/Vivacity_example.txt};
        \end{axis}
      \end{tikzpicture}
    \fi
  }
  \subfigure[Measured\label{sfig:Vivacity pixel 21_14}]{
    \iffinal
      \includegraphics{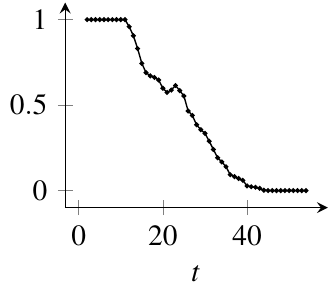}
    \else
      \begin{tikzpicture}
        \begin{axis}[%
          axis x line   = bottom,
          axis y line   = left,
          width         = 4.25cm,
          xlabel        = {$t$},
          enlargelimits = true,
        ]
          \addplot[sharp plot, mark=*, mark size=0.5pt] table[header=false] {Data/Vivacity/Vivacity_21_14.txt};
        \end{axis}
      \end{tikzpicture}
    \fi
  }
  \subfigure[Simulation\label{sfig:Vivacity pixel simulation}]{
    \iffinal
      \includegraphics{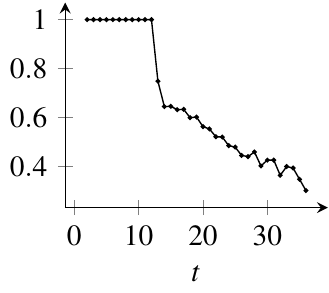}
    \else
      \begin{tikzpicture}
        \begin{axis}[%
          axis x line   = bottom,
          axis y line   = left,
          width         = 4.25cm,
          xlabel        = {$t$},
          enlargelimits = true,
        ]
          \addplot[sharp plot, mark=*, mark size=0.5pt] table[header=false] {Data/Vivacity/Vivacity_simulation.txt};
        \end{axis}
      \end{tikzpicture}
    \fi
  }
  \subfigureCaptionSkip
  \caption{%
    Vivacity curves~(pixel-based) for the three different data sets. At a glance, the curves permit
    comparing the dynamics of each process. 
    The sampling frequencies are different, necessitating the use of dynamic time warping.
  }
  \label{fig:Vivacity example}
\end{figure}
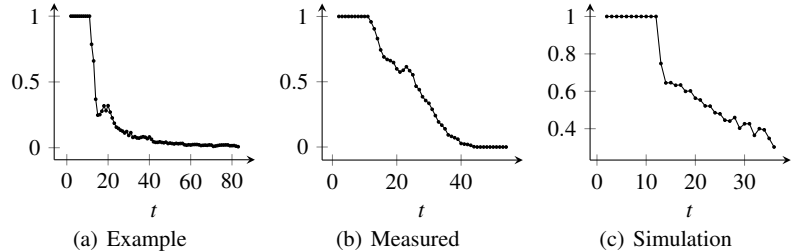

\runinhead{Quantifying dissimilarity}
%
To quickly quantify the dissimilarity between different curves, e.g., the vivacity curves that we
defined in Sec.~\ref{sec:Vivacity}, we can use \emph{dynamic time warping}~\cite{Berndt94},
a technique from dynamic programming that is able to compensate for different sampling frequencies
and different simulation lengths.
Fig.~\ref{fig:Vivacity example} depicts the vivacity curves of the data sets. We can see that the measured data in Fig.~\ref{sfig:Vivacity pixel 21_14} is characterized by a slower process in which
new mass is continuously being injected to the system. Hence, its vivacity does not decrease
steeply as that of the example data in Fig.~\ref{sfig:Vivacity pixel example}.
The vivacity curve for the simulation, shown in Fig.~\ref{sfig:Vivacity pixel simulation}, appears
to differ from the remaining curves. As a consequence, we can use these curves in visual comparison
tasks and distinguish between different~(measured) experiments and simulations.
The dynamic time warping distance helps quantify this assumption. We have $\dist(a,b) \approx 442$,
$\dist(a,c) \approx 1135$, and $\dist(b,c) \approx 173$. This indicates that the characteristics of
the simulation in Fig.~\ref{sfig:Vivacity pixel simulation} differ from those found in a real-world
viscous fingering process, shown in Fig.~\ref{sfig:Vivacity pixel example}, while being reasonably
close to another measured experiment, which is depicted by Fig.~\ref{sfig:Vivacity pixel 21_14}.
Vivacity curves may thus be used for parameter tuning of simulations in order to obtain better
approximations to measured data.

\section{Conclusion}

Driven by the need for a coherent analysis of time-varying skeletons, we developed different concepts
inspired by topological persistence in this paper.
We showed how to improve the consistency of tracking algorithms between consecutive time steps.
Moreover, we demonstrated the utility of our novel concepts for different purposes, including the
persistence-based filtering of skeletons, anomaly detection, and characterization of dynamic
processes.

Nonetheless, we envision numerous other avenues for future research.
For example, the propagation velocity of structures in the data may be of interest in many
applications.
We also plan to provide a detailed analysis of viscous fingering, including domain expert feedback,
and extend persistence to physical concepts within this context.
More generally, our novel persistence-inspired concepts can also be used in other domains, such as
the analysis of motion capture data~(which heavily relies on skeletonization techniques) or
time-varying point geometrical point clouds, for which novel skeletonization techniques
were recently developed~\cite{Kurlin15}.


\bibliographystyle{spmpsci}
\bibliography{TopoInVis2017_Viscous_Fingering}

\end{document}